\documentclass[reqno,a4paper]{amsart}
\usepackage{amssymb}
\usepackage{latexsym}
\usepackage{amsmath}

\input{diagrams}
\newcommand {\IN}{\mathbb{N}}                          
\newcommand {\IZ}{\mathbb{Z}}

\newcommand {\IR}{\mathbb{R}}                           
\newcommand {\IC}{\mathbb{C}}
\newcommand {\T}{\mathbb T}

\newcommand {\D}{\mathcal D}

\newcommand {\G}{\mathcal G}
\renewcommand {\H}{\mathcal H}

\newcommand {\K}{\mathcal K}
\renewcommand {\L}{\mathcal L}

\newcommand {\U}{\mathfrak U}

\newcommand {\g}{\mathfrak{g}}
\newcommand {\gc}{\mathfrak{g}_{_{\IC}}}
\renewcommand {\t}{\mathfrak{t}}
\newcommand {\tc}{\t_{\IC}}

\newcommand {\half}[1]{\frac{#1}{2}}

\newcommand {\rot}{\operatorname{Rot}(S^{1})}

\newtheorem {theorem}{Theorem}[subsection]
\newtheorem {proposition}[theorem]{Proposition}
\newtheorem {corollary}[theorem]{Corollary}
\newtheorem {lemma}[theorem]{Lemma}

\newcommand {\ssection}[1]{\section{#1}\setcounter{theorem}{0}\setcounter{equation}{0}}
\newcommand {\ssubsection}[1]{\subsection{#1}\hfill\vspace{0.5\baselineskip}}
\renewcommand {\thetheorem}
{\ifcase\arabic{section}\else\arabic{section}.\ifcase\arabic{subsection}\else\arabic{subsection}.\fi\fi\arabic{theorem}}

\numberwithin {equation}{subsection}
\renewcommand {\theequation}
{\ifcase\arabic{section}\else\arabic{section}.\ifcase\arabic{subsection}\else\arabic{subsection}.\fi\fi\arabic{equation}}

\renewcommand {\proof}{{\sc Proof.}\ }

\newcommand {\remark}{{\sc Remark.\ }}

\newcommand {\halmos}{$\Diamond$}

\newcommand {\lpol}{L^{\scriptscriptstyle{\operatorname{pol}}}}

\newcommand {\Hom}{\operatorname{Hom}}
\newcommand {\End}{\operatorname{End}}
\newcommand {\Ker}{\operatorname{Ker}}

\newcommand {\<}{\langle}
\renewcommand {\>}{\rangle}
\newcommand {\sqnm}[1]{\<#1,#1\>}

\newcommand {\pl}{\mathcal P_{\ell}}

\newcommand {\fin}{^{\scriptscriptstyle{\operatorname{fin}}}}
\newcommand {\Hfin}{\H\fin}
\newcommand {\hfin}{\H\fin}

\newcommand {\diff}{\operatorname{Diff}(S^{1})}

\newcommand {\der}{\left.\frac{d}{dt}\right|_{t=0}}

\renewcommand {\root}{\Lambda_{R}}
\newcommand {\weight}{\Lambda_{W}}
\newcommand {\domweight}{\weight^{+}}
\newcommand {\vvee}{^{\raisebox{.6pt}{$\scriptscriptstyle{\vee}$}}}
\newcommand {\coroot}{\Lambda_{R}\vvee}
\newcommand {\coweight}{\Lambda_{W}\vvee}
\newcommand {\lambdaz}{\Lambda_{Z}\vvee}

\newcommand {\primed}[1]{{#1}^{\prime}}
\newcommand {\cor}[1]{\alpha_{#1}^{\vee}}	
\newcommand {\cow}[1]{\lambda_{#1}^{\vee}}	
\newcommand {\cowp}[1]{\primed{\cow{#1}}}	

\newcommand {\lzg}{L_{Z}G}

\newcommand {\Ad}{\operatorname{Ad}}
\newcommand {\ad}{\operatorname{ad}}

\newcommand {\wt}{\widetilde}
\newcommand {\wh}{\widehat}

\newcommand {\al}{\mathcal A_{\ell}}

\newcommand {\deltabar}{\overline{\Delta}}

\newcommand {\SU}{\operatorname{SU}}
\newcommand {\Spin}{\operatorname{Spin}}
\renewcommand {\Sp}{\operatorname{Sp}}
\newcommand {\SO}{\operatorname{SO}}
\newcommand {\PSO}{\operatorname{PSO}}

\newcommand {\ie}{{\it i.e. }}

\newcommand {\cinfty}[1]{C^{\infty}(S^{1},#1)}
\newcommand {\LGl}{\wt{LG}^{\ell}}
\newcommand {\Lgl}{\wt{L\g}^{\ell}}

\renewcommand {\diff}{\operatorname{Diff}_{+}(S^{1})}
\newcommand {\vect}{\operatorname{Vect}(S^{1})}
\newcommand {\dtheta}{\frac{d\theta}{2\pi}}
\newcommand {\ddtheta}{\frac{d}{d\theta}}
\renewcommand {\der}[1]{\left.\frac{d}{d#1}\right|_{#1=0}}

\newcommand {\Int}{\int_{0}^{2\pi}}

\renewcommand {\P}{\mathcal P}

\newcommand {\lgz}{L(G/Z)}
\newcommand {\ellb}{\ell_{b}}
\newcommand {\ellf}{\ell_{f}}

\newcommand {\SpinDD}{\operatorname{Spin}_{4n}}
\newcommand {\twotwo}{\IZ_{2}\times\IZ_{2}}
\newcommand {\wl}[1]{\mathfrak{w}_{#1}}
\newcommand {\zigma}[1]{\mathfrak{s}_{#1}}

\newcommand {\slc}{\mathfrak{sl}_{2}(\IC)}
\newcommand {\soDD}{\mathfrak{so}_{4n}}
\newcommand {\soDDc}{\mathfrak{so}_{4n,\IC}}
\newcommand {\comm}[2]{#1#2#1^{-1}#2^{-1}}

\newcommand {\pback}{\pi^{*}U(\H)}
\newcommand {\pbackdue}{\rho^{*}U(\K)}
\newcommand {\pil}{\pi_{\lambda}}

\newcommand {\lyg}{L_{Y}G}
\newcommand {\ind}[2]{\operatorname{ind}_{#1}^{#2}}

\newcommand {\rotk}{\operatorname{Rot}^{k}(S^{1})}
\newcommand {\diffk}{\operatorname{Diff}^{k}_{+}(S^{1})}
\newcommand {\diffm}{\operatorname{Diff}^{m}_{+}(S^{1})}

\begin{document}

\begin{center}
{\tt Communications in Mathematical Physics 207 (1999), 307-339.}
\end{center}

\title
[Positive Energy Representations of Loop Groups]
{Positive Energy Representations of the Loop Groups of Non Simply Connected
Lie Groups}
\author{Valerio Toledano Laredo}
\address{
Institut de Math\'ematiques de Jussieu\\
UMR 7586\\
Case 191\\
Universit\'e Pierre et Marie Curie\\
4, Place Jussieu\\
F--75252 Paris Cedex 05}
\email{toledano@math.jussieu.fr}
\thanks{
Work supported by a TMR fellowship, contract no. FMBICT950083
and by a Non--Commutative Geometry Network grant, contract no.
ERB FMRXCT960073, both from the European Commission.}
\begin{abstract}
We classify and construct all irreducible positive energy
representations of the loop group of a compact, connected
and simple Lie group and show that they admit an intertwining
action of Diff$(S^{1})$.
\end{abstract}
\maketitle

\ssection{Introduction}
%=====================

Let $K$ be a compact, connected and simple Lie group and
$LK=C^{\infty}(S^{1},K)$ its loop group. We shall be
concerned with the study of positive energy representations
of $LK$, \ie projective unitary representations
\begin{equation}
\pi:LK\longrightarrow PU(\H)=U(\H)/\T
\end{equation}
extending to the semi--direct product $LK\rtimes\rot$ in
such a way that $\rot$ acts by non--negative characters
only and with finite--dimensional eigenspaces. In other
words,
\begin{equation}
\H=\bigoplus_{n\geq 0}\H(n)
\end{equation}
where
$\H(n)=\{\xi\in\H|\pi(R_{\theta})\xi=e^{in\theta}\xi\}$, the
subspace of energy $n$, supports a finite--dimensional
projective representation of $K$.\\

Positive energy representations are completely reducible and, when
$K$ is simply connected, have been classified by several authors
\cite{PS,Wa}. The irreducible ones are then uniquely determined
by their level $\ell\in\IN$ and lowest energy subspace $\H(0)$.
The former classifies the corresponding central extension of $LK$
or, equivalently, the cocycle associated to the infinitesimal
representation of its Lie algebra. The latter is an irreducible
$K$--module the highest weight $\lambda$ of which is bound by
the requirement that
\begin{equation}
\<\lambda,\theta\>\leq\ell
\end{equation}
where $\theta$ is the highest root of $K$ and $\<\cdot,\cdot\>$
is the basic inner product of $K$, \ie the multiple of the Killing
form such that $\<\theta,\theta\>=2$.\\

The aim of the present paper is to extend the above classification
to the case of groups which are not simply connected. Write $K=G/Z$
where $G$ is the universal covering group of $K$ and
$\pi_{1}(K)\cong Z\subseteq Z(G)$. It will be more convenient
to consider positive energy representations of the group of
discontinuous loops
\begin{equation}
\lzg=
\{\zeta\in C^{\infty}(\IR,G)|\zeta(x+2\pi)\zeta(x)^{-1}\in Z\}
\end{equation}
deferring to later the determination of those which factor
through $LK$.\\

Since $\lzg/LG\cong Z$, these may be studied with Mackey's
machine \cite{Ma1,Ma2}, paying however due care to the fact
that the representations in question are genuinely projective
and form a strict subclass of those of $\lzg$. With these
{\it provisos}, the analysis carries over essentially unchanged
and is dealt with in the following sections. We summarise them
below.\\

In section \ref{se:extensions of lzg}, we classify the central
extensions of $\lzg$ and $L(G/Z)$ by $\T$. We show in particular
the existence of an obstruction to the extension of the level
$\ell$ central extension of $LG$ to $\lzg$. This appears only
at odd $\ell$ and for some $G$, the complete list of which is
given in an appendix. It comprises $SU_{n}$ with $n$ even.
Thus, surprisingly perhaps, $L_{\IZ_{2}}\SU_{2}$, and {\it a
fortiori} $L\SO_{3}$ do not have any odd level positive energy
representations.
A further obstruction on the level appears when demanding that
a given central extension of $\lzg$ descend to $L(G/Z)$.
$\ell$ must then be a multiple of a given non--negative integer
$\ellb$, the basic level of $G/Z$, of which we compute the value
for all simple groups. It is $n$ for $\SU_{n}/\IZ_{n}$.\\

In section \ref{se:disc loops}, we show that the category $\pl$
of positive energy representations of $LG$ at a given level $\ell$
is closed under conjugation by elements of $\lzg$ and therefore
that $Z\cong\lzg/LG$ acts on the positive energy dual of $LG$.
We also compute the geometric counterpart of this action on the
alcove of $G$ which parametrises the irreducibles in $\pl$.
Aside from rendering the action of $Z$ more explicit, this shows
that it operates by automorphisms of the extended Dynkin diagram
of $G$.\\

In section \ref{se:Mackey obstruction}, we compute the Mackey
obstruction for a subgroup $Y\subseteq Z$ stabilising a given
positive energy representation $\H$ of $LG$. This vanishes for
most groups since $Y$ is cyclic unless $G/Y=\PSO_{4n}$ but,
somewhat surprisingly, doesn't in the latter case.\\

Section \ref{se:per of lzg} contains our main results. We construct
all irreducible positive energy representations of $\lzg$ and show
that they are classified by the central extension of $\lzg$ they
induce and their isomorphism class as $LG$--modules. Moreover,
we prove that they admit an intertwining action of $\diff$ and
identify it with the Segal--Sugawara representation obtained by
regarding them as positive energy $LG$--modules. Finally, in
section \ref{se:factoring} we determine those representations which
factor through $L(G/Z)$, thereby obtaining all positive energy
representations of the latter group. They are exactly those the
level of which is a multiple of the basic level of $G/Z$.\\

\remark
In physical terms, we classify in this paper all inequivalent
quantisations of the chiral Wess--Zumino--Witten model with
target group $G/Z$. Related results have been obtained in
the non--chiral case by Gepner--Witten \cite{GW},
Felder--Gaw\c{e}dzki--Kupiainen \cite{FGK1,FGK2} and
Gaberdiel \cite{Ga}.\\

{\bf Acknowledgements.}
I am grateful to A. Wassermann for suggesting the use of Mackey's
machine in the present context and for many useful conversations.
This paper was begun at the Department of Pure Mathematics and
Mathematical Statistics of the University of Cambridge and
completed at the Institut de Math\'ematiques de Jussieu of
the Universit\'e Pierre et Marie Curie, within the Alg\`ebres
d'Op\'erateurs et Repr\'esentations group. I wish to thank both
institutions for their kind hospitality and pleasant working
atmospheres.

\ssection{The coroot and coweight lattices of $G$}\label{se:lattices}
%================================================

We begin by gathering some elementary properties of the lattices canonically
associated to $G$. The present discussion follows \cite{GO}. Throughout this
paper, $G$ denotes a compact, connected and simply connected simple Lie group
with Lie algebra $\g$. Let $T\subset G$ be a maximal torus with Lie algebra
$\t\subset\g$. By the roots of $G$ we shall always mean its infinitesimal
roots, namely the set $R$ of linear forms $\alpha\in i\t^{*}=\Hom(\t,i\IR)$
such that the subspace
$\g_{\alpha}=\{x\in\gc|\thinspace[h,x]=\alpha(h)x\thickspace\forall h\in\tc\}$
is non--zero. Let $\Delta=\{\alpha_{1},\ldots,\alpha_{n}\}$ be a basis of
$R$ and $\theta$ the corresponding highest root.
The basic inner product $\<\cdot,\cdot\>$, \ie the unique multiple of the
Killing form such that $\<\theta,\theta\>=2$, is positive definite on $i\t$
and gives an identification $i\t^{*}\cong i\t$ of which we shall make implicit
use. The coroots  of $G$ are the elements of $i\t$ given by
$\alpha^{\vee}=\frac{2\alpha}{\<\alpha,\alpha\>}$. They form the dual root
system $R^{\vee}$.\\

The root and coroot lattices $\root\subset i\t^{*}$, $\coroot\subset i\t$ are
the lattices spanned by $R$ and $R^{\vee}$ respectively. They have $\IZ$--basis
given by $\Delta$ and $\Delta^{\vee}=\{\cor{1},\ldots,\cor{n}\}$. Since $\theta$
is a long root and there are at most two root lengths in $R$ with the ratio of
the squared length of a long root by that of a short one equal to 2 or 3,
rewriting $\alpha^{\vee}=\frac{\<\theta,\theta\>}{\<\alpha,\alpha\>}\alpha$
we see that $\coroot\subset\root$. Notice that 
$\<\alpha^{\vee},\alpha^{\vee}\>= \frac{4}{\<\alpha,\alpha\>}=
2\frac{\<\theta,\theta\>}{\<\alpha,\alpha\>}$
so that $\coroot$ is an even, and therefore integral lattice. The weight and
coweight lattices $\weight\subset i\t^{*}$,
$\coweight\subset i\t$ are the lattices dual to $\coroot$ and $\root$ respectively.
They have $\IZ$--basis given by the fundamental (co)weights $\lambda_{i}$,
$\cow{i}$ defined by
\begin{equation}
\<\lambda_{i},\cor{j}\>=\<\cow{i},\alpha_{j}\>=\delta_{ij}
\end{equation}
Clearly, $\coweight\subset\weight$. Moreover, by the integrality properties of
root systems, $\<\alpha,\beta^{\vee}\>\in\IZ$ for any root $\alpha$ and coroot
$\beta^{\vee}$ so that $\root\subset\weight$ and, dually,
$\coroot\subset\coweight$. Graphically,

\begin{equation}\label{eq:summa}
\begin{diagram}[height=2em,width=3em]
\root   &\subset&\weight        &\subset        &i\t^{*}\\
\cup    &       &\cup           &               &       \\
\coroot &\subset&\coweight      &\subset        &i\t    \\
\end{diagram}
\end{equation}

Let $Z(G)$ be the centre of $G$ and $\wh{Z(G)}=\Hom(Z(G),\T)$ its Pontryagin
dual. The following is well--known

\begin{lemma}\label{centre}\hfill
\begin{enumerate}
\item The map $e(h)=\exp_{T}(-2\pi ih)$ induces an isomorphism
$\coweight/\coroot\cong Z(G)$.
\item The pairing $\mu(\exp_{T}(h))=e^{\<\mu,h\>}$ induces an
isomorphism $\weight/\root\cong\wh{Z(G)}$.
\end{enumerate}
\end{lemma}
% \proof
% (i) Since $G$ is connected, $T$ is maximal abelian and therefore $Z(G)\subset T$.
% It follows that $Z(G)\cong e^{-1}(Z(G))/\Ker e$ where the integral lattice
% $\Ker e\cong\Hom(\T,T)$ is equal to $\coroot$ since $G$ is simply--connected
% \cite[thm. 5.47] {Ad}.
% To show that $e(h)\in Z(G)$ iff $h\in\coweight$, we use $Z(G)=\Ker(\Ad_{G})$ and
% the fact that if $0\neq x_{\alpha}\in\g_{\alpha}$,
% then
% $\Ad(e(h))x_{\alpha}=\exp(-2\pi i\ad(h))x_{\alpha}=e^{-2\pi i\alpha(h)}x_{\alpha}$
% is equal to $x_{\alpha}$ iff $\alpha(h)\in\IZ$.\\
% (ii) The map $\weight/\root\rightarrow\wh{\coweight/\coroot}$,
% $\mu\rightarrow e^{-2\pi i\<\mu,\cdot\>}$ is readily seen to be an isomorphism
% and coincides with the given pairing under the identification
% $Z(G)\cong\coweight/\coroot$
% \halmos\\

\remark When $G$ is simply--laced, \ie with all roots of equal length,
the basic inner product identifies roots and coroots and the vertical
inclusions in \eqref{eq:summa} are equalities. Moreover, lemma \ref{centre}
yields a canonical isomorphism $\wh{Z(G)}\cong Z(G)$.\\

The Weyl group $W$ of $G$ is the finite group generated in $\End(i\t^{*})$
by the orthogonal reflections $\sigma_{\alpha}$ corresponding to the roots
$\alpha\in R$. Since
\begin{equation}
 \sigma_{\alpha}(\mu)=
 \mu-2\frac{\<\mu,\alpha\>}{\<\alpha,\alpha\>}\alpha=
 \mu-\<\mu,\alpha^{\vee}\>\alpha=
 \mu-\<\mu,\alpha\>\alpha^{\vee}
\end{equation}
the action of $W$ preserves $\coroot$--cosets in $\coweight$ and $\root$--cosets
in $\weight$. Call $\mu\in\weight$ (resp. $\mu\in\coweight$) minimal if it
is of minimal length in its $\root$ (resp. $\coroot$)--coset. The following
gives a characterisation of minimal (co)weights.

\begin{proposition}\label{pr:minimal orbit}
There is, in each $\coweight/\coroot$--coset (resp. $\weight/\root$--coset) a unique
$W$--orbit of elements of minimal length. These may equivalently be characterised as
those $\lambda$ such that
\begin{equation}\label{minimal}
\<\lambda,\alpha\>\in\{0,\pm 1\}\qquad
\text{(resp. $\<\lambda,\alpha^{\vee}\>\in\{0,\pm 1\}$)}
\end{equation} for any root $\alpha$ (resp. coroot $\alpha^{\vee}$).
\end{proposition}
\proof
It is sufficient to consider the case of $\coweight/\coroot$ since $\root,\weight$
are the coroot and coweight lattices of the dual root system $R^{\vee}$.
Let $\mu\in\coweight$ be of minimal length in its $\coroot$--coset. Then, for any
root $\beta$ and corresponding coroot $\beta^{\vee}=\frac{2\beta}{\<\beta,\beta\>}$,
we have $\|\mu\pm\beta^{\vee}\|^{2}\geq\|\mu\|^{2}$ and, expanding
$|\<\mu,\beta\>|\leq 1$. Assume that $\lambda\in\coweight$ satisfies \eqref{minimal}
and $\nu=\lambda$ mod $\coroot$ is of minimal length in its coset. We claim that
$w\lambda=\nu$ for an appropriate $w\in W$. To see this, write
$\nu=\lambda+\beta^{\vee}_{1}+\cdots+\beta^{\vee}_{r}$ where the $\beta^{\vee}_{i}$
are (possibly repeated) coroots. Clearly, one cannot have
$\<\lambda,\beta_{i}\>\geq 0$ for all $i$ otherwise
\begin{equation}
\<\nu,\nu\>=
\<\lambda,\lambda\>+\<\sum\beta^{\vee}_{i},\sum\beta^{\vee}_{i}\> +
2\<\lambda,\sum\beta^{\vee}_{i}\> >
\<\lambda,\lambda\>
\end{equation}
in contradiction with the minimality of $\nu$. Thus, by \eqref{minimal} there
exists an $i\in\{1,\ldots,r\}$ such that $\<\lambda,\beta_{i}\>=-1$ and therefore
$\lambda_{1}:=\sigma_{\beta_{i}}\lambda=\lambda+\beta^{\vee}_{i}$. Moreover,
$\lambda_{1}$ satisfies \eqref{minimal} since $W$ permutes the roots and preserves
$\<\cdot,\cdot\>$. We may therefore iterate the above step to find a permutation
$\tau$ of $\{1,\ldots,r\}$ such that
\begin{equation}
\lambda_{i}:=
\lambda+\beta^{\vee}_{\tau(1)}+\cdots+\beta^{\vee}_{\tau(i)}=
\sigma_{\beta_{\tau(i)}}\cdots\sigma_{\beta_{\tau(1)}}\lambda
\end{equation}
In particular, $\lambda_{r}=\nu$ and therefore $\nu\in W\lambda$ whence
$\|\lambda\|=\|\nu\|$ \halmos\\

Recall that a weight $\mu\in\weight$ is dominant if it lies in the cone
\begin{equation}
 \domweight=
 \{\nu\in\weight|\thinspace\<\nu,\cor{i}\>\geq 0
 \thickspace\forall\cor{i}\in\Delta^{\vee}\}=
 \bigoplus\lambda_{i}\cdot\IN
\end{equation}
Since $\domweight$ is a fundamental domain for the action of $W$ on $\weight$,
lemma \ref{centre} and proposition \ref{pr:minimal orbit} establish a bijective
correspondence between elements in $\wh{Z(G)}$ and minimal dominant weights.
Dually, the elements of $Z(G)$ correspond to the minimal dominant coweights,
\ie those $\mu\in(\coweight)^{+}=\bigoplus_{i}\cow{i}\cdot\IN$ of minimal
length in their $\coroot$--coset. The following gives another characterisation
of minimal dominant coweights.

\begin{lemma}\label{special roots}
The non--zero minimal dominant coweights are exactly the fundamental
co\-weights corresponding to special roots, \ie those $\alpha_{i}\in\Delta$
bearing the coefficient 1 in the expansion
\begin{equation}\label{highest root}
\theta=\sum m_{i}\alpha_{i}
\end{equation}
\end{lemma}
\proof
By proposition \ref{pr:minimal orbit}, $\mu\in(\coweight)^{+}$ is
minimal iff $\<\mu,\theta\>\leq 1$. Indeed, for any positive root
$\alpha$, we get
$0\leq\<\mu,\alpha\>\leq\<\mu,\theta\>-\<\mu,\theta-\alpha\>
  \leq\<\mu,\theta\>$.
Since $\<\mu,\theta\>=0$ implies $\mu=0$, the non--zero minimal
dominant coweights are those $\mu\in(\coweight)^{+}$ such that
$\<\mu,\theta\>=1$. Writing $\mu=\sum_{i}k_{i}\cow{i}$, $k_{i}\geq 0$
and using \eqref{highest root}, we find $\<\mu,\theta\>=\sum k_{i}m_{i}$.
Since $\theta-\alpha_{i}$ is a sum of positive roots, $m_{i}\geq 1$
for any $i$ and result follows \halmos

\ssection{Central extensions of $\lzg$} \label{se:extensions of lzg}
%=====================================

This section is devoted to the study of the central extensions by
$\T$ of the group of discontinuous loops
\begin{equation}
\lzg=\{\zeta\in C^{\infty}(\IR,G)|\medspace
       \zeta(x+2\pi)\zeta(x)^{-1}\in Z\}
\end{equation}
corresponding to a subgroup $Z\subseteq Z(G)$. These are uniquely
determined by their restrictions to $LG=\cinfty{G}$ and to
$\Hom(\T,T/Z)$, the integral lattice of $G/Z$. The former are
classified by their level $\ell\in\IZ$ \cite{PS} and the latter
by their commutator map, a $\T$--valued, skew--symmetric bilinear
form $\omega$ on $\Hom(\T,T/Z)$. We shall prove below that $\ell$
and $\omega$ are bound by the requirement that
\begin{equation}\label{eq:intro}
\omega(\lambda,\mu)=(-1)^{\ell\<\lambda,\mu\>}
\end{equation}
whenever $\lambda$ lies in the coroot lattice $\Hom(\T,T)$ and
therefore that central extensions of $LG$ do not necessarily
extend to $\lzg$ since a suitable $\omega$ satisfying
\eqref{eq:intro} for a given $\ell$ need not exist. In particular,
$L_{\IZ_{2}}SU_{2}$, and more generally $L_{\IZ_{2n}}SU_{2n}$, do
not possess central extensions of odd level.
For compatible $\ell$ and $\omega$, we construct the corresponding
central extension of $\lzg$ and show that the action of $\diff$
on $\lzg$ lifts uniquely to it. The classification of central
extensions of $L(G/Z)$ follows easily from this and is
described at the end of this section.

\ssubsection{Central extensions of $\mathbf{LG}$}
\label{ss:frechet nonsense}
%------------------------------------------------

We begin by reviewing the construction of central extensions of
$LG$, and more generally of a connected and simply--connected
Fr\'echet Lie group $\G$, following chapter 4 of \cite{PS}.
All central extensions considered in this section are understood
to be smooth and have $\T$ as their extending group. Let $\L$ be
the Lie algebra of $\G$ and $\beta$ a two--cocycle on $\L$, \ie
a continuous, skew--symmetric, bilinear map
$\beta:\L\times\L\rightarrow\IR$ satisfying
\begin{equation}
\beta([X,Y],Z)+\beta([Y,Z],X)+\beta([Z,X],Y)=0
\end{equation}
$\beta$ may be regarded as a right--invariant, closed two--form
on $\G$ and we assume that $(2\pi)^{-1}\beta$ is integral, \ie
such that its integral over any two--cycle in $\G$ is an integer.
Then, there exists a unique central extension
\begin{equation}
1\rightarrow\T\rightarrow\wt\G\xrightarrow{\pi}\G\rightarrow 1
\end{equation}
the Lie algebra of which is $\wt\L=\L\oplus i\IR$ with bracket
\begin{equation}\label{connection}
[X\oplus it,Y\oplus is]=[X,Y]\oplus i\beta(X,Y)
\end{equation}

$\wt\G$ may be constructed using the following path group description.
Assume $\wt\G$ exists and regard it as a principal $\T$--bundle over
$\G$ with connection given by the splitting $\wt\L=\L\oplus i\IR$. In
other words, the horizontal subspace at $\wt g\in\wt\G$ is $\L\wt g$.
The pull--back of $\wt\G$ to the space
\begin{equation}
\P\G=\{p:I\rightarrow \G|\thinspace p(0)=1\}
\end{equation}
of piece--wise smooth paths via the end--point fibration
$\P\G\xrightarrow{e}\G$
is topologically trivial, the identification of the fibre at the
constant path 1 with that at $p$ being simply given by parallel
transport along $p$. Explicitly, if $X=\dot pp^{-1}:I\rightarrow\L$
is the right logarithmic derivative of $p$, the identification maps
$z\in\T=e^{*}\pi^{-1}(1)$ to the end point of the path $\wt p$ in
$\wt\G$ obtained by solving $\dot{\wt p}=X\wt p$, $\wt p(0)=z$.
If $p$ is closed, and therefore contractible in $\G$, the corresponding
identification is simply multiplication by the holonomy
$e^{i\int_{\sigma}\beta}$ where $\sigma$ is any two--cycle in $\G$ with
boundary $p$.\\

The concatenation of pointed paths defined by
\begin{equation}
p\vee q(t)=
\left\{\begin{array}{rcl}
q(2t)&\text{if}&0\leq t\leq\half{1}\\
p(2t-1)q(1)&\text{if}&\half{1}\leq t\leq 1
\end{array}\right.
\end{equation}
induces a monoidal structure on $\P\G$ which, combined with
the group law on $\wt\G$ makes $e^{*}\wt\G$ a monoid. The
crucial feature of the corresponding multiplication law is
that it becomes the canonical one when transported to
$\P\G\times\T\cong e^{*}\wt\G$, a direct consequence
of the $\wt\G$--invariance of the connection on $\wt\G$. It
follows that, as a group, $\wt\G$ may be described, or indeed
{\it defined} as the quotient of $\P\G\times\T$ with law
$(p,z)\star(q,w)=(p\vee q,zw)$ by the equivalence relation
\begin{equation}
(p,z)\sim(q,w)
\quad\Longleftrightarrow\quad
p(1)=q(1)
\thickspace\thickspace\text{and}\thickspace\thickspace
e^{i\int_{\sigma}\beta}=w\overline{z}
\end{equation}
where $\sigma$ is a two--cycle with boundary $p\vee\check{q}$
and $\check{q}(t)=q(1-t)q(1)^{-1}$.

\begin{lemma}\label{le:lift}
An automorphism $A$ of $\G$ lifts to $\wt\G$ if, and only if it
leaves the cohomology class of $\beta$ invariant, \ie iff there
exists a linear map $F:\L\rightarrow\IR$ such that for any $X,Y
\in\L$
\begin{equation}\label{coboundary}
\beta(AX,AY)=\beta(X,Y)+F([X,Y])
\end{equation}
The lift is then unique up to multiplication by a character of
$\G$ and is given infinitesimally by
\begin{equation}\label{inf lift}
\wt{A}(X\oplus it)=AX\oplus i(F(X)+t)
\end{equation}
and in the path group description of $\wt\G$ by
\begin{equation}\label{phase correction}
\wt{A} (p,z)=(Ap,ze^{i\int_{p}F})
\end{equation}
where $F$ is regarded as a right--invariant one--form, so that
$\int_{p}F=\int_{0}^{1}F(\dot pp^{-1})$.
\end{lemma}
\proof
The necessity of \eqref{coboundary} is straightforward. Indeed,
a lift $\wt{A}$ acts on $\wt\L$ by
$\wt{A}(X\oplus it)=AX\oplus i(G(X)+t)$ for some linear map
$G:\L\rightarrow\IR$. Requiring that
$\wt{A}[X\oplus i,Y\oplus i]=[\wt{A}(X\oplus i),\wt{A}(Y\oplus i)]$
and expanding both members yields
\begin{equation}
A[X,Y]\oplus i(G([X,Y])+\beta(X,Y))=[AX,AY]\oplus i\beta(AX,AY)
\end{equation}
and therefore \eqref{coboundary}. Conversely, \eqref{phase correction}
is a well--defined lift of A. Indeed, when regarded as an identity
between right--invariant forms in $\G$, \eqref{coboundary} reads
$A^{*}\beta=\beta-dF$. It follows that if $(p,z)\sim (q,w)$ and
$\sigma$ is a two--cycle in $\G$ with $\partial\sigma=p\vee\check{q}$,
then, $\partial\medspace A\sigma=Ap\vee\check{(Aq)}$ and
\begin{equation}
e^{i\int_{A\sigma}\beta}=e^{i\int_{\sigma}A^{*}\beta}=
e^{i\int_{\sigma}\beta}e^{-i\int_{p}F}e^{i\int_{q}F}=
we^{i\int_{q}F}\overline{ze^{i\int_{p}F}}
\end{equation}
so that $(Ap,ze^{i\int_{p}F})\sim (Aq,we^{i\int_{q}F})$.
The uniqueness of $\wt{A}$ is clear for if $\wt{A}_{i}$, $i=1,2$
are two lifts of $A$, then $\wt{A}_{2}\wt{A}_{1}^{-1}$ is
a lift of the identity and fixes $\T$ so that it is given
by $\chi\circ\pi$ for some $\chi\in\Hom(\G,\T)$ \halmos\\

\remark The phase factor in \eqref{phase correction} may be
derived from \eqref{inf lift} as follows. For any $p\in\P\G$,
denote by $\wt p$ its unique horizontal lift through $1\in\wt\G$
so that $\dot{\wt p}=\dot pp^{-1}\wt p$ and $\wt p(0)=1$.
Then $Q(t)=\wt A\wt p(t)$ solves
\begin{equation}\label{eq:inho}
\dot{Q}=
\wt{A}(\dot pp^{-1})Q=
\dot{(Ap)}{(Ap)}^{-1}Q+iF(\dot{p}p^{-1})Q
\end{equation}
Set $Q(t)=\phi(t)\wt{Ap}(t)$ where $\phi(t)\in\T$, then
\eqref{eq:inho} reduces to $\dot\phi=iF(\dot pp^{-1})\phi$
and therefore $\phi(t)=e^{i\int_{0}^{t}F(\dot pp^{-1})d\tau}$.
Conversely, \eqref{inf lift} may be obtained from
\eqref{phase correction} by taking $p$ as the path
$s\rightarrow\exp_{\G}(stX)$ and differentiating at $t=0$.\\

Let now $\G=LG=\cinfty{G}$ with Lie algebra $L\g=\cinfty{\g}$.
The basic inner product $\<\cdot,\cdot\>$ on $\g$ determines a
right--invariant, closed two--form on $LG$ given by
\begin{equation}\label{eq:basic cocycle}
B(X,Y)=\int_{0}^{2\pi}\<X,\dot Y\>\frac{d\theta}{2\pi}
\end{equation}
and such that $(2\pi)^{-1}B$ is integral. For any $\ell\in\IZ$,
denote by $\LGl$ and $\Lgl$ the central extension of $LG$ corresponding
to $\ell B$ and its Lie algebra. Then, any central extension of
$LG$ is isomorphic to some $\LGl$ for a uniquely determined
$\ell\in\IZ$ called its level \cite[thm. 4.4.1]{PS}.\\

Since $B$ is invariant under the action of the group $\diff$
of orientation--preserving diffeomorphisms of $S^{1}$ given
by $\phi\gamma=\gamma\circ\phi^{-1}$, and $\Hom(LG,\T)=\{1\}$
\cite[prop. 3.4.1]{PS}, this action lifts uniquely to any
central extension of $LG$.
Similarly, the action by conjugation of $\lzg$ on $LG$ lifts
to any $\LGl$. Indeed, for $\zeta\in\lzg$ we have
\begin{equation}\label{eq:zeta on B}
\begin{split}
B(\zeta X\zeta^{-1},\zeta Y\zeta^{-1})
&=
\Int\<\zeta X\zeta^{-1},\zeta\dot Y\zeta^{-1}\>\dtheta+
\Int\<\zeta X\zeta^{-1},\zeta[\zeta^{-1}\dot\zeta,Y]\zeta^{-1}\>\dtheta\\
&=
B(X,Y)-\Int\<\zeta^{-1}\dot\zeta,[X,Y]\>\dtheta
\end{split}
\end{equation}
where we used the Ad--invariance of $\<\cdot,\cdot\>$ and
the fact that
$\dot{\zeta^{-1}}\zeta+\zeta^{-1}\dot{\zeta}=
 \dot{(\zeta^{-1}\zeta)}=0$.
Therefore, by lemma \ref{le:lift}, on $\Lgl$
\begin{equation}\label{eq:zeta on lie}
\wt{\Ad(\zeta)}\medspace X\oplus it =
\zeta X\zeta^{-1}\oplus i
\left(t-\ell\Int\<\zeta^{-1}\dot\zeta,X\>\dtheta\right)
\end{equation}
In particular, the adjoint action of $\LGl$ factors through
$LG\subset\lzg$ and, by the uniqueness of lifts, is given by
\eqref{eq:zeta on lie}.

\ssubsection{The compatibility requirement}
\label{ss:compatibility}
%------------------------------------------

By lemma \ref{centre}, the subgroup $Z\subseteq Z(G)$ is isomorphic
to $\lambdaz/\coroot$ where $\coroot\subset\lambdaz\subset\coweight$
is the integral lattice of $G/Z$, \ie $\lambdaz\cong\Hom(\T,T/Z)$.
We will regard $\lambdaz$ as a subgroup of $\lzg$ by associating
to $\mu\in\lambdaz$ the discontinuous loop
$\zeta_{\mu}(\theta)=\exp_{T}(-i\theta\mu)$.
Since any character of $\coroot$ extends to $\lambdaz$, the
connecting homomorphism in the five term sequence
\begin{equation}
\Hom(\lambdaz,\T)\rightarrow\Hom(\coroot,\T)
\rightarrow H^{2}(Z,\T)\rightarrow H^{2}(\lambdaz,\T)
\end{equation}
is the zero map. Using the five term sequence for the inclusion
$LG\subset\lzg$ and the fact that $\Hom(LG,\T)=\{1\}$
\cite[Prop. 3.4.1]{PS},
we therefore obtain the following commutative diagram with exact row

\begin{equation}\label{comm exact}
\begin{diagram}[height=1.8em,width=3em]
 &    &0                 &     &              &    &            \\
 &    &\dTo              &     &              &    &            \\
0&\rTo&H^{2}(Z,\T)       &\rTo &H^{2}(\lzg,\T)&\rTo&H^{2}(LG,\T)\\
 &    &\dTo              &\ldTo&              &    &            \\
 &    &H^{2}(\lambdaz,\T)&     &              &    &            \\
\end{diagram}
\end{equation}

which shows that a central extension of $\lzg$ by $\T$ is entirely
determined by its restrictions to $LG$ and to $\lambdaz$. The former
is classified by its level $\ell$ and the latter by its commutator
map $\omega$ defined by
\begin{equation}\label{eq:def of omega}
\omega(\lambda,\mu)=
\wt\zeta_{\lambda}\wt\zeta_{\mu}
\wt\zeta_{\lambda}^{-1}\wt\zeta_{\mu}^{-1}
\end{equation}
where $\wt\zeta_{\lambda},\wt\zeta_{\mu}\in\wt\lzg$ are arbitrary
lifts of $\zeta_{\lambda}, \zeta_{\mu}$. $\omega$ is a skew--symmetric,
$\T$--valued, $\IZ$--bilinear form on $\lambdaz$. Since the identification
$\lambdaz\cong\Hom(\T,T/Z)$ maps $\coroot$ to $\Hom(\T,T)\subset LG$,
$\omega$ is bound by the requirement that
$\omega(\alpha,\beta)=(-1)^{\ell\<\alpha,\beta\>}$
whenever $\alpha,\beta\in\coroot$ \cite[prop. 4.8.1]{PS}. We shall
presently establish that $\omega$ is constrained by a more astringent
identity, the proof of which gives an alternative derivation of
proposition 4.8.1 of \cite{PS}.

\begin{theorem}\label{th:constraint}
Let $\wt{\lzg}$ be a central extension of $\lzg$ by $\T$, the
restrictions to $LG$ and $\lambdaz$ of which have level $\ell$
and commutator map $\omega$ respectively. Then, for any
$\lambda\in\coroot$ and $\mu\in\lambdaz$ 
\begin{equation}\label{commutator}
\omega(\mu,\lambda)=(-1)^{\ell\<\mu,\lambda\>}
\end{equation}
\end{theorem}

Theorem \ref{th:constraint} is an immediate corollary of the
following

\begin{proposition}\label{pr:konstraint}
For any $\mu\in\coweight$, denote by $\wt A_{\mu}$ the unique lift
of the conjugation action of $\zeta_{\mu}$ on $LG$ to $\LGl$. Then,
for any $\lambda\in\coroot$ and lift $\wt\zeta_{\lambda}\in\LGl$ of
$\zeta_{\lambda}$,
\begin{equation} \label{kommutator}
\wt A_{\mu}(\wt\zeta_{\lambda})=
(-1)^{\ell\<\mu,\lambda\>}\wt\zeta_{\lambda}
\end{equation}
\end{proposition}

{\sc Proof (of theorem \ref{th:constraint}).}
Let $\wt\zeta_{\lambda},\wt\zeta_{\mu}\in\wt\lzg$ be lifts
of $\zeta_{\lambda}$ and $\wt\zeta_{\mu}\in\wt\lzg$ respectively.
By lemma \ref{le:lift}, $\Ad(\wt\zeta_{\mu})=\wt A_{\mu}$ as
automorphisms of $\LGl\cong\left.\wt\lzg\right|_{LG}$ since
both are lifts of $\Ad(\zeta_{\mu})$ and $\Hom(LG,\T)=\{1\}$.
Thus, by \eqref{eq:def of omega} and \eqref{kommutator}
\begin{equation}
\omega(\mu,\lambda)=
\wt A_{\mu}(\wt\zeta_{\lambda})\wt\zeta_{\lambda}^{-1}=
(-1)^{\ell\<\mu,\lambda\>}
\end{equation}
\halmos\\

{\sc Proof (of proposition \ref{pr:konstraint}).}
Since $\zeta_{\mu}\zeta_{\lambda}\zeta_{\mu}^{-1}=\zeta_{\lambda}$
in $LG$, the left hand--side of \eqref{kommutator} is equal to
$\omega(\mu,\lambda)\wt\zeta_{\lambda}$, where $\omega(\mu,\lambda)\in\T$
is independent of the choice of the lift $\wt\zeta_{\lambda}$ and
is bilinear in $\mu,\lambda$. Moreover, if $\mu\in\coroot$, lemma
\ref{le:lift} implies that $\wt A_{\mu}=\Ad(\wt\zeta_{\mu})$ and
$\omega$ is therefore skew--symmetric when restricted to
$\coroot\times\coroot$. We begin by establishing that
\begin{equation}\label{eq:to show}
\omega(\mu,\lambda)=(-1)^{\ell\<\mu,\lambda\>}
\end{equation}
when $\lambda=\alpha^{\vee}$ is the coroot corresponding to
a positive root $\alpha$ and $\mu\in\coweight$ is such that
$\<\alpha,\mu\>\in\{0,1\}$. The loop
$\zeta_{\alpha^{\vee}}(\theta)=\exp_{T}(-i\theta\alpha^{\vee})$
may then be written as a product of two exponentials in $LG$
\cite[4.8.1]{PS}, namely
\begin{equation}\label{product of exp}
\zeta_{\alpha^{\vee}}=
\exp_{LG}\Bigl(-\half{\pi}(e_{\alpha}(0)-f_{\alpha}(0) )\Bigr)
\exp_{LG}\Bigl( \half{\pi}(e_{\alpha}(1)-f_{\alpha}(-1))\Bigr)
\end{equation}
Here, using standard notation, $e_{\alpha},f_{\alpha}$ and
$h_{\alpha}=\alpha^{\vee}$ span the
$\mathfrak{sl}_{2}(\IC)$--subalgebra of $\gc$ corresponding
to $\alpha$ and, for any $x\in\gc$ and $n\in\IN$,
$x(n)=x\otimes e^{in\theta}\in L\gc$.
To see that \eqref{product of exp} holds, consider the homomorphism
$\sigma_{\alpha}:SU_{2}\rightarrow G$ mapping the standard basis of
$\mathfrak{sl}_{2}(\IC)$ given by
\begin{xalignat}{3}
e&=\begin{pmatrix}0&1\\0&0 \end{pmatrix}&
f&=\begin{pmatrix}0&0\\1&0 \end{pmatrix}&
h&=\begin{pmatrix}1&0\\0&-1\end{pmatrix}
\end{xalignat}
to $\{e_{\alpha},f_{\alpha},h_{\alpha}\}$. This induces a
homomorphism $LSU_{2}\rightarrow LG$ sending
\begin{equation}
\theta\longrightarrow
\begin{pmatrix}e^{-i\theta}&0\\0&e^{i\theta}\end{pmatrix}
\end{equation}
to $\zeta_{\alpha^{\vee}}$ and \eqref{product of exp} reduces to a
simple matrix check.\\

If $h\in\t$, then $[h,e_{\alpha}]=\<h,\alpha\>e_{\alpha}$ whence
$\Ad(\exp_{T}(h))e_{\alpha}=\exp(\ad(h))e_{\alpha}=
 e^{\<h,\alpha\>}e_{\alpha}$.
Therefore, since $\zeta_{\mu}(\theta)=\exp_{T}(-i\theta\mu)$, we
have
\begin{equation}
 \zeta_{\mu}e_{\alpha}(n)\zeta_{\mu}^{-1}(\theta)=
 e_{\alpha}\otimes e^{i\theta(n-\<\alpha,\mu\>)}=
 e_{\alpha}(n-\<\alpha,\mu\>)(\theta)
\end{equation} so that, in $L\g$
\begin{align}
\zeta_{\mu}e_{\alpha}(n)\zeta_{\mu}^{-1}&=e_{\alpha}(n-\<\alpha,\mu\>)\\
\zeta_{\mu}f_{\alpha}(n)\zeta_{\mu}^{-1}&=f_{\alpha}(n+\<\alpha,\mu\>)
\end{align}
Since $\zeta_{\mu}^{-1}\dot\zeta_{\mu}=-i\mu\in\t$ and this subspace
is orthogonal to $\IC e_{\alpha}\oplus\IC f_{\alpha}$ with respect 
to the Killing form, no correction term arises from \eqref{eq:zeta on lie}
and the same holds in $\Lgl$. It follows that \eqref{eq:to show} holds
if $\<\alpha,\mu\>=0$. If, on the other hand $\<\alpha,\mu\>=1$, then
\begin{equation}\label{temp}
\begin{split}
\wt A_{\mu}(\wt\zeta_{\alpha^{\vee}})\wt\zeta_{\alpha^{\vee}}^{-1}
=&
\exp_{\LGl}\Bigl(-\half{\pi}(e_{\alpha}(-1)-f_{\alpha}(1) )\Bigr)
\exp_{\LGl}\Bigl( \half{\pi}(e_{\alpha}(0) -f_{\alpha}(0) )\Bigr)\\
\cdot&
\exp_{\LGl}\Bigl(-\half{\pi}(e_{\alpha}(1) -f_{\alpha}(-1))\Bigr)
\exp_{\LGl}\Bigl( \half{\pi}(e_{\alpha}(0) -f_{\alpha}(0) )\Bigr)\\
=&
\exp_{\LGl}\Bigl(-\half{\pi}(e_{\alpha}(-1)-f_{\alpha}(1) )\Bigr)\\
\cdot&\Ad\left(
\exp_{\LGl}\Bigl( \half{\pi}(e_{\alpha}(0) -f_{\alpha}(0 ))\Bigr)\right)
\exp_{\LGl}\Bigl(-\half{\pi}(e_{\alpha}(1) -f_{\alpha}(-1))\Bigr)\\
\cdot&
\exp_{\LGl}\Bigl(       \pi (e_{\alpha}(0) -f_{\alpha}(0) )\Bigr)
\end{split}
\end{equation}
As is readily checked using $\sigma_{\alpha}$, we have
\begin{equation}
\Ad\left(
\exp_{LG}\Bigl( \half{\pi}(e_{\alpha}(0) -f_{\alpha}(0) )\Bigr)\right)
(e_{\alpha}(1)-f_{\alpha}(-1))
%\exp_{LG}\Bigl( \half{\pi}(e_{\alpha}(0) -f_{\alpha}(0) )\Bigr)^{-1}
=e_{\alpha}(-1)-f_{\alpha}(1)
\end{equation}
Moreover, since we are conjugating by a constant loop, no correction
term arises from \eqref{eq:zeta on lie} and \eqref{temp} is therefore
equal to
\begin{equation}\label{intermediate}
\exp_{\LGl}\Bigl(-\pi(e_{\alpha}(-1)-f_{\alpha}(1))\Bigr)
\exp_{\LGl}\Bigl( \pi(e_{\alpha}(0) -f_{\alpha}(0))\Bigr)
\end{equation}

To proceed, we seek to diagonalise the above elements. This is best
done in $LSU_{2}$ using the identity
\begin{equation}
e_{\alpha}(-m)-f_{\alpha}(m)=V(m)ih_{\alpha}(0)V(m)^{*}
\end{equation}
where $V(m)\in LSU_{2}$ is given by
$\theta\longrightarrow\frac{1}{\sqrt{2}}
 \begin{pmatrix}1&ie^{-im\theta}\\ie^{im\theta}&1\end{pmatrix}$.
Since
\begin{equation}
V^{-1}(m)\dot V(m)=
\half{m}\Bigl(ih_{\alpha}(0)+e_{\alpha}(-m)-f_{\alpha}(m)\Bigr)
\end{equation}
we have
\begin{equation}
\int_{0}^{2\pi}\<V^{-1}(m)\dot V(m),ih_{\alpha}(0)\>\frac{d\theta}{2\pi}=
-\half{m}\|h_{\alpha}\|^{2}=-m\frac{2}{\<\alpha,\alpha\>}
\end{equation}
and therefore, using \eqref{eq:zeta on lie}, \eqref{intermediate}
is equal to
\begin{equation}
e^{-i\pi\ell\frac{2}{\<\alpha,\alpha\>}}
\wt{V(1)}\exp_{\LGl}(-i\pi h_{\alpha}(0))\wt{V(1)}^{-1}
\wt{V(0)}\exp_{\LGl}( i\pi h_{\alpha}(0))\wt{V(0)}^{-1}
\end{equation}
where $\wt{V(0)},\wt{V(1)}$ are arbitrary lifts of $V(0),V(1)$ in $\LGl$.
Since $\exp_{SU_{2}}(-i\pi h_{\alpha}(0))=-1$ lies in the centre of any
central extension of $LSU_{2}$, the above is equal to
$(-1)^{\ell\frac{2}{\<\alpha,\alpha\>}}=(-1)^{\ell\<\alpha^{\vee},\mu\>}$
and \eqref{eq:to show} holds if $\<\alpha,\mu\>=1$.\\

Let now $\lambda=\alpha^{\vee}$ and $\mu=\beta^{\vee}$ be coroots. Then,
either $|\<\alpha,\beta^{\vee}\>|\leq 1$ or $|\<\alpha^{\vee},\beta\>|\leq 1$
\cite[Table 1, p.45]{Hu}. Using the bilinearity and skew--symmetry of both
sides of \eqref{eq:to show}, we may assume, up to a permutation and a
sign change, that $\alpha$ is positive and that
$\<\alpha,\beta^{\vee}\>\in\{0,1\}$ so that, by the computation above,
\eqref{eq:to show} holds whenever $\lambda$ and $\mu$ lie in the coroot
lattice. To complete the proof, it is sufficent to check \eqref{eq:to show}
when $\lambda=\alpha^{\vee}$ is a positive coroot and $\mu$ varies in a set
of representatives of $\coroot$--cosets in $\coweight$. A convenient choice
is given by the minimal dominant coweights. If $\mu$ is one such then,
by proposition \ref{pr:minimal orbit}, $\<\mu,\alpha\>\in\{0,1\}$ and
therefore \eqref{eq:to show} holds by our previous computation \halmos

\ssubsection{Construction of central extensions of $\mathbf{\lzg}$}
\label{ss:lzg tilde}
%------------------------------------------------------------------

Define the level $\ell$ and commutator map $\omega$ of a central
extension $\wt\lzg$ of $\lzg$ by restriction to $LG$ and $\lambdaz$
respectively. By theorem \ref{th:constraint}, no such $\wt\lzg$ exists
unless $\ell$ and $\omega$ are compatible, \ie satisfy \eqref{commutator}.
In particular, $L_{\IZ_{2}}SU_{2}$, and {\it a fortiori} $LSO_{3}$,
do not possess any central extensions of odd level since in this
case $\coroot=\alpha\IZ$ and $\coweight=\half{\alpha}\IZ$ with
$\<\alpha,\alpha\>=2$. On the other hand, \eqref{commutator}
requires $\omega(\alpha,\half{\alpha})=-1$ in contradiction
with the skew--symmetry of $\omega$.\\

Let now $\ell\in\IZ$ and $\omega$ be a skew--symmetric, $\T$--valued
bilinear form on $\lambdaz$. Then,

\begin{proposition}\label{pr:ext of lzg}
There exists a (necessarily unique) central extension $\wt\lzg$ of
$\lzg$ of level $\ell$ and commutator map $\omega$ if, and only if
\begin{equation}\label{eq:ns cond}
\omega(\mu,\lambda)=(-1)^{\ell\<\mu,\lambda\>}
\end{equation}
whenever $\lambda\in\coroot$.
\end{proposition}
\proof
The necessity of \eqref{eq:ns cond} is the contents of
theorem \ref{th:constraint} and the uniqueness of $\wt
\lzg$ that of \eqref{comm exact}.
Let $\wt{LG}$ and $\wt{\lambdaz}$ be the central extensions
of $LG$ and $\lambdaz$ with level $\ell$ and commutator map
$\omega$ respectively. Following \cite[Prop. 4.6.9]{PS}, we
shall construct $\wt\lzg$ as a quotient of
$\wt{LG}\rtimes\wt{\lambdaz}$
\footnote{The proof of proposition 4.6.9 of \cite{PS} is
sligthly erroneous in that it does not assume that $\omega$
is compatible with $\ell$. When that is not the case, the
group $N$ defined by \eqref{enn} is not normal as can be
seen from equation \eqref{normality two}.}.
Lift the conjugation action of $\lambdaz\subset\lzg$ on $LG$
to $\wt{LG}$ by using lemma \ref{le:lift} and denote the
corresponding automorphisms of $\wt{LG}$ by $\wt{A}_{\mu}$,
$\mu\in\lambdaz$. Form the semi--direct product
$\wt{LG}\rtimes\wt{\lambdaz}$ where the action of
$\wt{\lambdaz}$ factors trough $\lambdaz$. By theorem
\ref{th:constraint} and the compatibility of $\ell$ and $\omega$,
$\wt{LG}$ and $\wt{\lambdaz}$ restrict to isomorphic central
extensions $\wt{\coroot}$ of $\coroot$. We may therefore consider
the subgroup
\begin{equation}\label{enn}
N=\{(\wt\zeta_{\alpha},{\wt\zeta_{\alpha}}^{-1})\}
\subset\wt{LG}\rtimes\wt{\lambdaz}
\end{equation}
where $\wt\zeta_{\alpha}$ varies in $\wt{\coroot}$. We claim
that $N$ is normal. By lemma \ref{le:lift}, for any $\alpha\in\coroot$,
$\wt{A}_{\alpha}=\Ad(\wt\zeta_{\alpha})$ since both are lifts of
$\Ad(\zeta_{\alpha})$. Therefore, for any $\wt\gamma\in\wt{LG}$,
\begin{equation}
(\wt\gamma,1)
(\wt\zeta_{\alpha},{\wt\zeta_{\alpha}}^{-1})(\wt\gamma^{-1},1)=
(\wt\gamma,1)
(\wt\zeta_{\alpha}\wt{A}_{-\alpha}(\wt\gamma^{-1}),{\wt\zeta_{\alpha}}^{-1})=
(\wt\zeta_{\alpha},{\wt\zeta_{\alpha}}^{-1})
\end{equation}
Moreover, by proposition \ref{pr:konstraint} and \eqref{eq:ns cond}
\begin{equation}\label{normality two}
(1,\wt\zeta_{\mu})
(\wt\zeta_{\alpha},{\wt\zeta_{\alpha}}^{-1})(1,{\wt\zeta_{\mu}}^{-1})=
(\wt{A}_{\mu}(\wt\zeta_{\alpha}),
[\wt\zeta_{\mu},{\wt\zeta_{\alpha}}^{-1}]{\wt\zeta_{\alpha}}^{-1})=
((-1)^{\ell\<\mu,\alpha\>}\wt\zeta_{\alpha},
 (-1)^{-\ell\<\mu,\alpha\>}{\wt\zeta_{\alpha}}^{-1})
\end{equation}
Thus, $N$ is normal, and the quotient $\wt{LG}\rtimes\wt{\lambdaz}/N$
is the required central extension of $\lzg$ \halmos\\

\remark Theorem \ref{th:constraint} and proposition \ref{pr:ext of lzg}
prove the exactness of
\begin{equation}
1\rightarrow H^{2}(Z,\T)\rightarrow H^{2}(\lzg,\T)\rightarrow H^{2}(LG,\T)
 \rightarrow\IZ/\ell_{f}\IZ\rightarrow 0
\end{equation}
where $\ell_{f}$ is 1 if $\lambdaz$ possesses a commutator map
satisfying \eqref{commutator} with $\ell=1$ and $\ell_{f}=2$
otherwise. We call $\ell_{f}$ the fundamental level of
$G/Z$. The fundamental levels of all compact simple groups are
given in \S \ref{ss:ellb and ellf}.

\ssubsection{Automorphic action of $\mathbf{\diff}$ on $\mathbf{\wt\lzg}$}
\label{ss:action of diff}
%-------------------------------------------------------------------------

Let $\diff$ be the group of orientation--preserving diffeomorphisms
of $S^{1}$ and $\D$ its universal covering group. $\D$ may be
realised as the subgroup of diffeomorphisms $\phi$ of $\IR$
such that $\phi(x+2\pi)=\phi(x)+2\pi$ and $\diff\cong\D/(T_{2\pi})$
where $T_{y}$ is translation by $y$. $\D$ acts automorphically
on $\lzg$ by
\begin{equation}\label{eq:diff on lzg}
\phi\medspace\zeta=\zeta_{\phi^{-1}}=\zeta\circ\phi^{-1}
\end{equation}
and this action factors to one of $\D/(T_{2\pi k})$ where
$k$ is the order of the largest cyclic subgroup of $Z$.
The Lie algebra of $\D$ is the Lie algebra $\vect$ of all
smooth vector fields on $S^{1}$ with bracket
\footnote{\noindent
The bracket \eqref{eq:vect bracket} is the Lie--theoretic
bracket on $\vect$ satisfying
\begin{equation*}
[X,Y]=\der{t}\exp(tX)Y\exp(-tX)
\end{equation*}
and is the opposite of the differential geometric one
defined by the action of $\vect$ on $C^{\infty}(S^{1})$.}
\begin{equation}\label{eq:vect bracket}
[f\ddtheta,g\ddtheta]=(\dot{f}g-f\dot{g})\ddtheta
\end{equation}
If $\xi=f\ddtheta\in\vect$, the action of $\xi$ on $L\g$
corresponding to \eqref{eq:diff on lzg} is simply
\begin{equation}
\xi X=-f\dot{X}
\end{equation}
so that the Lie algebra of $\D\rtimes\lzg$ is
$\vect\rtimes L\g$ with bracket
\begin{equation}\label{eq:semidir bracket}
[X\oplus f\ddtheta,Y\oplus g\ddtheta]=
([X,Y]-f\dot Y+g\dot X)\oplus (\dot{f}g-f\dot{g})\ddtheta
\end{equation}

\begin{proposition}\label{pr:action of diff}
Let $k$ be the order of the largest cyclic subgroup of $Z$.
Then, the action of the universal $k$--covering of $\diff$
on $\lzg$ lifts uniquely to any central extension $\wt\lzg$.
\end{proposition}
\proof
The uniqueness is easily settled for two lifts necessarily
differ by some
\begin{equation}
\chi\in\Hom(\D,\Hom(\lzg,\T))\cong\Hom(\D,\wh{Z})=1
\end{equation}
where the first isomorphism follows from $\Hom(LG,\T)=1$ and
the second by connectedness of $\D$.
Let $\wt{LG},\wt{\lambdaz}$ be the restrictions of $\wt\lzg$
to $LG$ and $\lambdaz$ respectively and $\ell,\omega$ the
corresponding level and commutator map. We shall describe,
as in the proof of proposition \ref{pr:ext of lzg}, $\wt\lzg$
as a quotient of $\wt{\lambdaz}\ltimes\wt{LG}$.
Since $\lzg$ isn't connected, the action of $\D$ on it
cannot be lifted to $\wt\lzg$ by using lemma \ref{le:lift}.
Following the proof of \cite[prop. 4.7.1]{PS}, we shall
regard it instead as one of $\lzg$ on $\wt{LG}\rtimes\D$
in the following way.\\

Consider first the connected component of the identity of
$\wt\lzg$, \ie $\wt{LG}$ and form the semi--direct product
$\wt{LG}\rtimes\D$ where $\D$ acts on $\wt{LG}$ as in \S
\ref{ss:frechet nonsense}. Since $\D$ is contractible, $\wt{LG}\rtimes\D$
may equivalently be described as the central extension of $LG\rtimes\D$
corresponding to the Lie algebra cocycle
\begin{equation}
\beta(X\oplus f\ddtheta,Y\oplus g\ddtheta)=
\ell B(X,Y)=\ell\Int\<X,\dot Y\>\dtheta
\end{equation}

We claim that $\lzg$ acts on $LG\rtimes\D$ and $\wt{LG}\rtimes\D$.
The first action simply stems from the fact that $LG\rtimes\D$ is
a normal subgroup of $\lzg\rtimes\D$ and is given explicitly by
\begin{equation}\label{eq:lzg on LGD}
\zeta(\gamma,\phi)=
(\zeta\gamma\zeta^{-1}\zeta\zeta^{-1}_{\phi^{-1}},\phi)
\end{equation}
and infinitesimally by
\begin{equation}\label{eq:lzg on LgD}
\zeta(X\oplus f\ddtheta)=
(\zeta X\zeta^{-1}+f\dot{\zeta}\zeta^{-1})\oplus f\ddtheta
\end{equation}

To see that this action lifts to $\wt{LG}\rtimes\D$, we compute
\begin{equation}
\beta(\zeta(X\oplus f\ddtheta),\zeta(Y\oplus g\ddtheta))=
\ell B(\zeta X\zeta^{-1}+f\dot\zeta\zeta^{-1},
       \zeta Y\zeta^{-1}+g\dot\zeta\zeta^{-1})
\end{equation}
By \eqref{eq:zeta on B},
\begin{equation}
B(\zeta X\zeta^{-1},\zeta Y\zeta^{-1})=
B(X,Y)-\Int\<[X,Y],\zeta^{-1}\dot\zeta\>\dtheta
\end{equation}

On the other hand,

\begin{equation}
\begin{split}
B(f\dot\zeta\zeta^{-1},\zeta Y\zeta^{-1})+
B(\zeta X\zeta^{-1},g\dot\zeta\zeta^{-1})
&=\Int f\<\dot\zeta\zeta^{-1},\zeta\dot{Y}\zeta^{-1}\>
 +     f\<\dot\zeta\zeta^{-1},\zeta[\zeta^{-1}\dot\zeta,Y]\zeta^{-1}\>\dtheta\\
&-\Int g\<\dot\zeta\zeta^{-1},\zeta\dot{X}\zeta^{-1}\>
 -     g\<\dot\zeta\zeta^{-1},\zeta[\zeta^{-1}\dot\zeta,X]\zeta^{-1}\>\dtheta\\
%&=\Int f\<\zeta^{-1}\dot\zeta,\dot{Y}\>\dtheta
% +\Int f\<\zeta^{-1}\dot\zeta,[\zeta^{-1}\dot\zeta,Y]\>\dtheta\\
%&-\Int g\<\zeta^{-1}\dot\zeta,\dot{X}\>\dtheta
% -\Int g\<\zeta^{-1}\dot\zeta,[\zeta^{-1}\dot\zeta,X]\>\dtheta\\
&=\Int\<\zeta^{-1}\dot\zeta,f\dot{Y}-g\dot{X}\>\dtheta
\end{split}
\end{equation}

Finally, anti--symmetrising, we find
\begin{equation}
B(f\dot\zeta\zeta^{-1},g\dot\zeta\zeta^{-1})=
\half{1}\Int(f\dot{g}-\dot{f}g)
\<\dot\zeta\zeta^{-1},\dot\zeta\zeta^{-1}\>\dtheta
\end{equation}

Thus
\begin{equation}
\beta(\zeta(X\oplus f\ddtheta),\zeta(Y\oplus g\ddtheta))=
\beta(X\oplus f\ddtheta,Y\oplus g\ddtheta)-
\ell F([X\oplus f\ddtheta,Y\oplus g\ddtheta])
\end{equation}
where $F:L\g\rtimes\vect\rightarrow\IR$ is given by
\begin{equation}
F(X\oplus f\ddtheta)=
\Int\<X,\zeta^{-1}\dot\zeta\>\dtheta+
\half{1}\Int f\<\zeta^{-1}\dot\zeta,\zeta^{-1\dot\zeta}\>\dtheta
\end{equation}

Since $\D$ is perfect \cite{Ep}, it follows from lemma \ref{le:lift} that the
action of $\lzg$ on $LG\rtimes\D$ lifts uniquely to $\wt{LG}\rtimes\D$ and is
given by

\begin{equation}\label{eq:segal formula +}
\begin{split}
\zeta(X\oplus f\ddtheta\oplus it)
&=
(\zeta X\zeta^{-1}+f\dot{\zeta}\zeta^{-1})\oplus f\ddtheta\\
&\oplus i\left(t-
\ell\Int\<X,\zeta^{-1}\dot\zeta\>\dtheta-
\half{\ell}\Int f\<\zeta^{-1}\dot\zeta,\zeta^{-1}\dot\zeta\>\dtheta
\right)
\end{split}
\end{equation}

Consider now the semi--direct product
$\wt{\lambdaz}\ltimes(\wt{LG}\rtimes\D)$ where the action of
$\wt{\lambdaz}$ factors through $\lambdaz$ and the subgroup 
\begin{equation}
N=\{(\wt{\zeta_{\alpha}},\wt{\zeta_{\alpha}}^{-1},1)\}
\subset\wt{\lambdaz}\ltimes(\wt{LG}\rtimes\D)
\end{equation}
where $\wt{\zeta_{\alpha}}$ varies in
$\wt{\coroot}=\left.\wt{\lzg}\right|_{\coroot}$. $N$ lies in
the centraliser of $\wt{LG}\rtimes\D$ since, by uniqueness, 
$\Ad((\wt{\zeta_{\alpha}},1,1))=\Ad((1,\wt{\zeta_{\alpha}},1))$
on $\wt{LG}\rtimes\D$ as both automorphisms are lifts of
$\Ad(\zeta_{\lambda})$. It follows that the quotient
$\wt{\lambdaz}\ltimes\wt{LG}/N\cong\wt\lzg$ is acted upon by
$\D$.\\

To conclude, we need only show that translations by multiples
of $2\pi k$ act trivially on $\wt\lzg$, where $k$ is the order
of the largest cyclic subgroup of $Z$. It is sufficient to check
this on a representative of each connected component of $\wt\lzg$
since, by uniqueness, $T_{2\pi}\wt{\gamma}=\wt\gamma$ for any
$\wt\gamma\in\wt{LG}$. Let
$\wt{\zeta_{\lambda}}$ be a lift of the discontinuous loop
$\zeta_{\lambda}(\theta)=\exp_{T}(-i\lambda\theta)$,
$\lambda\in\lambdaz$. By \eqref{eq:segal formula +},
\begin{equation}\label{eq:lzg on rot}
\wt{\zeta_{\lambda}}\ddtheta\wt{\zeta_{\lambda}}^{-1}=
-i\lambda\oplus\ddtheta\oplus i\half{\ell}\<\lambda,\lambda\>
\end{equation}
and therefore, since $T_{y}=\exp_{\D}(y\ddtheta)$
\begin{equation}\label{eq:integral on rot}
\wt{\zeta_{\lambda}}T_{2\pi k}\wt{\zeta_{\lambda}}^{-1}=
e^{\pi ik\ell\<\lambda,\lambda\>}\exp_{T}(-2\pi ik\lambda)T_{2\pi k}=
(-1)^{\ell\<k\lambda,\lambda\>}T_{2\pi k}
\end{equation}
Notice that $k\lambda\in\coroot$ since its image in $Z$ is 1.
Thus, by the skew--symmetry of $\omega$ and its compatibility
with $\ell$, we have 
\begin{equation}
1=\omega(k\lambda,\lambda)=(-1)^{k\ell\<\lambda,\lambda\>}
\end{equation}
whence
\begin{equation}
\wt{\zeta_{\lambda}}T_{2\pi k}\wt{\zeta_{\lambda}}^{-1}=
T_{2\pi k}
\end{equation}
as claimed \halmos\\

Let us record the following by--product of the proof of proposition
\ref{pr:action of diff} since it extends formula 4.9.4 of \cite{PS}

\begin{corollary}\label{cor:segal formulae}
The action of $\lzg$ on the Lie algebra of $\wt{LG}\rtimes\diff$,
where $\wt{LG}$ is the central extension of $LG$ of level $\ell$,
is given by
\begin{equation}
\begin{split}
\zeta\left(X\oplus f\ddtheta\oplus it\right)
&=
\left(\zeta X\zeta^{-1}+f\dot{\zeta}\zeta^{-1}\right)\oplus f\ddtheta\\
&\oplus
i\left(
t
-\ell\Int\<X,\zeta^{-1}\dot{\zeta}\>\dtheta
-\half{\ell}\Int f\<\zeta^{-1}\dot{\zeta},\zeta^{-1}\dot{\zeta}\>\dtheta
\right)
\end{split}
\end{equation}
\end{corollary}

\ssubsection{Central extensions of $\mathbf{\lgz}$}
\label{ss:ext of lgz}
%--------------------------------------------------

We now classify the central extensions of $\lgz$. The five
term sequence corresponding to
\begin{equation}
1\rightarrow Z\rightarrow\lzg\xrightarrow{\pi}\lgz\rightarrow 1
\end{equation}
and the fact that $\Hom(LG,\T)=1$ yield the exactness of
\begin{equation}
1\rightarrow\Hom(Z,\T)\rightarrow H^{2}(\lgz,\T)
 \xrightarrow{\pi^{*}}H^{2}(\lzg,\T)
\end{equation}
The image of $\pi^{*}$ is easily described.
Let the basic level $\ellb$ of $G/Z$ be the smallest integer
$\ell$ such that the restriction of $\ell\<\cdot,\cdot\>$ to
$\lambdaz$ is integral, \ie such that
\begin{equation}\label{eq:basic req}
\ell\<\cow{i},\cow{j}\>\in\IZ
\end{equation}
for all fundamental coweights $\cow{i},\cow{j}$ lying in
$\lambdaz$. Then,

\begin{proposition}\label{pr:ext of lgz}
A central extension of $\lzg$ is the pull--back of one of
$\lgz$ only if its level $\ell$ is a multiple of the basic
level of $G/Z$. Conversely, if $\ellb|\ell$, the subgroup
$Z\subset\wt\lzg$ corresponding to the canonical embedding
$G\hookrightarrow\wt\lzg$ is central and
\begin{equation}
\wt\lzg\cong\pi^{*}(\wt\lzg/Z)
\end{equation}
\end{proposition}
{\sc Proof}
\footnote{
The 'only if' implication of proposition \ref{pr:ext of lgz}
is essentially the contents of lemma 4.6.3 of \cite{PS}.}.
As readily verified, a central extension $\wt\lzg$ of
$\lzg$ is the pull--back of one of $\lgz$ only if its
restriction to $Z$ lies in its centre.
Conversely, since $G$ is simple and simply--connected,
the restriction of $\wt\lzg$ to $G$, and therefore to
$Z$, is canonically split. If
$s:Z\longrightarrow\wt Z=\left.\wt\lzg\right|_{Z}$ is
the corresponding section and $\wt Z$ is central, then
$\wt\lzg/s(Z)$ is a central extension of $\lgz$ which
pulls back to $\wt\lzg$.
We therefore need to determine those $\wt\lzg$
for which $\wt Z$ is a central subgroup. Notice
first that $\wt Z$ lies in the centre of
$\wt{LG}=\left.\wt\lzg\right|_{LG}$.
Indeed, for any $z\in Z$, $\gamma\in LG$
and lifts $\wt z,\wt\gamma\in\wt\lzg$,
\begin{equation}
\wt\gamma\wt z{\wt\gamma}^{-1}{\wt z}^{-1}=
\chi(\gamma,z)
\end{equation}
where $\chi(\gamma,z)\in\T$ is independent of the lifts
and multiplicative in each variable. Since
$\Hom(LG,\T)=\{1\}$ however, $\chi\equiv 1$.
Thus, we need only check that $\wt Z$ commutes with the
lifts of the discontinuous loops
$\zeta_{\lambda}(\theta)=\exp_{T}(-i\lambda\theta)$,
$\lambda\in\lambdaz$. For any $h\in\t$, we have by
\eqref{eq:zeta on lie},
\begin{equation}
{\wt\zeta_{\lambda}}h{\wt\zeta_{\lambda}}^{-1}=
h+\ell\<h,\lambda\>
\end{equation}
whence
\begin{equation}
{\wt\zeta_{\lambda}}\wt z{\wt\zeta_{\lambda}}^{-1}=
{\wt z}e^{-2\pi i\ell\<\mu,\lambda\>}
\end{equation}
where $\mu\in\lambdaz$ is such that $\exp_{T}(-2\pi i\mu)=z$,
and it follows that $\wt Z$ is central iff $\ell$ is a multiple
of $\ellb$ \halmos\\

Proposition \ref{pr:ext of lgz} shows that
\begin{equation}
H^{2}(\lgz,\T)\cong
\Ker\ell\oplus\Hom(Z,\T)\subseteq
H^{2}(\lzg,\T)\oplus\Hom(Z,\T)
\end{equation}
where $\ell$ is the map giving the residue mod $\ellb$ of
the level of a central extension of $\lzg$. The isomorphism
is simply given by associating to a central extension $\wt\lzg$
of level $\ell\in\ellb\IZ$ and $\chi\in\wh{Z}$ the central
extension
\begin{equation}
\wt\lzg/(z\cdot\chi(z))_{z\in Z}
\end{equation}
The list of basic levels for all compact, connected and simple
Lie groups is given in \S \ref{ss:ellb and ellf}.\\

\remark
If the central extension $\wt\lzg$ has level $\ell\in\ellb\IZ$,
the action of $\D$ on $\wt\lzg$ clearly descends to $\wt{L(G/Z)}=\wt\lzg/Z$.
Surprisingly perhaps, it does not then necessarily factor to one
of $\diff$. Indeed, \eqref{eq:lzg on rot} yields
\begin{equation}
T_{2\pi}\wt\zeta_{\lambda}T_{2\pi}^{-1}=
\wt\zeta_{\lambda}\exp_{T}(2\pi i\lambda)
e^{2\pi i\ell\half{\sqnm{\lambda}}}
\end{equation}
which equals 1 in $\wt\lzg/Z$ if, and only if
$\ell\sqnm{\lambda}\in 2\IZ$ for any $\lambda$ \ie iff
$\lambdaz$, endowed with $\<\cdot,\cdot\>$, is an even
lattice.\\

\remark The basic level of $G/Z$ is a multiple of the fundamental one
for if $\ellb|\ell$, the form
\begin{equation}
\omega(\lambda,\mu)=
(-1)^{\ell\<\lambda,\mu\>+\ell^{2}\<\lambda,\lambda\>\<\mu,\mu\>}
\end{equation}
is a commutator map on $\lambdaz$ satisfying the hypothesis
of proposition \ref{pr:ext of lzg}. In particular, $\lzg$
possesses a canonical central extension at level $\ellb$.

\ssubsection{Appendix : fundamental and basic levels of simple Lie groups}
\label{ss:ellb and ellf}
%-------------------------------------------------------------------------

Let $Z\subseteq Z(G)$ and $\coroot\subseteq\lambdaz\subseteq\coweight$
be the fundamental group and integral lattice of $G/Z$.

\begin{lemma}\label{le:cyclic}
If $Z\cong\lambdaz/\coroot$ is cyclic of order $k$, then $G/Z$ has
fundamental level 1 if and only if $k\<\lambda,\lambda\>\in 2\IZ$
where $\lambda\in\lambdaz$ is a generator.
\end{lemma}
\proof
If $\omega$ is a commutator map on $\lambdaz$ satisfying
\begin{equation}\label{eq:one commutator}
\omega(\alpha,\mu)=(-1)^{\<\alpha,\mu\>}
\end{equation}
whenever $\alpha\in\coroot$, then, by skew--symmetry,
$1=\omega(k\lambda,\lambda)=(-1)^{k\<\lambda,\lambda\>}$ since
$k\lambda\in\coroot$. Conversely, if $k\<\lambda,\lambda\>\in 2\IZ$,
the form
$\wt\omega(\alpha\oplus a\lambda,\beta\oplus b\lambda)=
 (-1)^{\<\alpha,\beta\>+\<b\alpha+a\beta,\lambda\>}$
on $\coroot\oplus\IZ\lambda$ descends to one on
$\coroot\oplus\IZ\lambda/-k\lambda\oplus k\lambda\cong\lambdaz$
satisfying \eqref{eq:one commutator} \halmos

\begin{proposition}\label{pr:ellb & ellf}
The following is the list of fundamental and basic levels
$\ellf,\ellb$ for all compact, connected and simple Lie
groups with universal cover $G$ and fundamental group
$Z\neq\{1\}$.

\begin{center}
\renewcommand {\arraystretch}{1.35}
\begin{tabular}{|| l l | l | l | l | l | l ||}\hline
{\bf G} &               &{\bf Z(G)}&{\bf Z}  &{\bf G/Z}  &$\mathbf{\ellf}$&$\mathbf{\ellb}$ \\\hline
%SU$_{n}$& $n\geq 2$     &$\IZ_{n}$ &$\IZ_{k}$&           &1 if $\frac{n(n-1)}{k}\in 2\IZ$&smallest $\ell$ with\\\cline{6-6}
%        &               &          &         &           &2 otherwise     &$\frac{n(n-1)}{k^{2}}\ell\in\IZ$\\\hline
SU$_{n}$& $n\geq 2$     &$\IZ_{n}$ &$\IZ_{k}$&           &1 for $n$ odd        &smallest $\ell$ with            \\
        &               &          &         &           &or $\frac{n}{k}$ even&$\frac{n(n-1)}{k^{2}}\ell\in\IZ$\\\cline{6-6}
        &               &          &         &           &2 otherwise          &                                \\\hline
Spin$_{2n+1}$& $n\geq 2$&$\IZ_{2}$ &$\IZ_{2}$&SO$_{2n+1}$&1               &1              \\\hline
Sp$_{n}$& $n\geq 1$     &$\IZ_{2}$ &$\IZ_{2}$&           &1 for $n$ even  &1 for $n$ even \\\cline{6-7}
        &               &          &         &           &2 for $n$ odd   &2 for $n$ odd  \\\hline
Spin$_{4m}$& $m\geq 2$  &$\IZ_{2}\times\IZ_{2}$&$\IZ_{2}^{0}$       &SO$_{4m}$ &1              &1 \\\cline{4-7}
        &               &                    &$\IZ_{2}^{\pm}$       &          &1 for $m$ even &1 for $m$ even \\\cline{6-7}
        &               &                    &                      &          &2 for $m$ odd  &2 for $m$ odd  \\\cline{4-7}
        &               &                    &$\IZ_{2}\times\IZ_{2}$&PSO$_{4m}$&1 for $m$ even &2              \\\cline{6-6}
        &               &                    &                      &          &2 for $m$ odd  &               \\\hline
Spin$_{4m+2}$& $m\geq 1$&$\IZ_{4}$ &$\IZ_{2}$&SO$_{4m+2}$&1               &1 \\\cline{4-7}
        &               &          &$\IZ_{4}$&PSO$_{4m+2}$&2              &4 \\\hline
E$_{6}$ &               &$\IZ_{3}$ &$\IZ_{3}$&           &1               &3 \\\hline
E$_{7}$ &               &$\IZ_{2}$ &$\IZ_{2}$&           &2               &2 \\\hline
\end{tabular}
\end{center}
\end{proposition}
\proof
We proceed by enumeration according to the Lie--theoretic type
of $G$, using the tables \cite[planches I--IX]{Bou} and lemmas
\ref{le:cyclic} and \ref{special roots}. For $G$ simply--laced,
we identify the coroot and coweight lattices with the root and
weight lattices respectively. In what follows, $\theta_{i}$,
$i=1\ldots n$ and $\<\cdot,\cdot\>$ are the standard basis and
inner product in $\IR^{n}$. Unless otherwise indicated, the
basic inner product is the standard one. 

${\bf SU_{n}}$, ${\bf n\geq 2}$\\
$\SU_{n}$ is simply--laced and the quotient $\weight/\root\cong\IZ_{n}$ is
generated by $\cow{1}=\theta_{1}-\frac{1}{n}(\sum_{i=1}^{n}\theta_{i})$
corresponding to the special root $\alpha_{1}=\theta_{1}-\theta_{2}$.
For any $k|n$, the subgroup of $Z(\SU_{n})$ isomorphic to $\IZ_{k}$
is generated by $\frac{n}{k}\cow{1}$ and
$\sqnm{\frac{n}{k}\cow{1}}=\frac{n(n-1)}{k^{2}}$.

{\bf{Spin$_{\bf 2n+1}$}}, ${\bf n\geq 2}$\\
$Z(\Spin_{2n+1})\cong\IZ_{2}$ is generated by the coweight $\cow{1}=
\theta_{1}$ corresponding to the unique special root $\alpha_{1}=
\theta_{1}-\theta_{2}$. Since $\sqnm{\cow{1}}=1$, $\ellb=\ellf=1$.

${\bf Sp_{n}}$, ${\bf n\geq 1}$\\
$\Sp_{1}$ is the group of unit quaternions and is therefore isomorphic
to $SU_{2}$. For $n\geq 2$, $Z(\Sp_{n})\cong\IZ_{2}$ is generated by
$\cow{n}=\theta_{1}+\cdots+\theta_{n}$ corresponding to the unique
special root $\alpha_{n}=2\theta_{n}$. Since the basic inner product
is half the standard one on $\IR^{n}$, we have $\sqnm{\cow{n}}=\half{n}$
whence $\ell_{b}=\ellf=1$ for $n$ even and $2$ for $n$ odd. This is
consistent with the isomorphism $\Sp_{2}\cong\Spin_{5}$.

{\bf{Spin$_{\bf 2n}$}, ${\bf n\geq 3}$}\\
Spin$_{2n}$ is simply--laced with minimal dominant coweights
$\cow{1}=\theta_{1}$,
$\cow{n-1}=\half{1}(\theta_{1}+\cdots+\theta_{n-1}-\theta_{n})$
and $\cow{n}=\half{1}(\theta_{1}+\cdots+\theta_{n})$ corresponding
to the special roots $\alpha_{1}=\theta_{1}-\theta_{2}$,
$\alpha_{n-1}=\theta_{n-1}-\theta_{n}$ and
$\alpha_{n}=\theta_{n-1}+\theta_{n}$. $2\cow{1}=0$ mod $\root$
and $\sqnm{\cow{1}}=1$ so that the corresponding quotient
$\Spin_{2n}/\IZ_{2}\cong\SO_{2n}$ has $\ellb=\ellf=1$.
We must now distinguish two cases :

{\it $n$ odd}. Then $2\cow{n-1}=2\cow{n}=\cow{1}$ mod $\root$
and $Z(\Spin_{2n})\cong\IZ_{4}$ with $\cow{n-1},\cow{n}$ of
order 4. Since $\sqnm{\cow{n}}=\frac{n}{4}$, we get $\ellb=4$
and $\ellf=2$ for $\Spin_{2n}/\IZ_{4}$. This is in agreement
with the isomorphism $\Spin_{6}\cong SU_{4}$.

{\it $n$ even}. Then $2\cow{n-1}=2\cow{n}=0$ mod $\root$ and
$Z(\Spin_{2n})\cong\IZ_{2}\times\IZ_{2}$. Since
$\sqnm{\cow{n-1}}=\sqnm{\cow{n}}=\frac{n}{4}$ we get
$\ellb=\ellf=1$ or $\ellb=\ellf=2$ for the quotients
$\Spin_{2n}/\IZ_{2}^{\pm}$ corresponding to $\cow{n-1}$ and
$\cow{n}$ according to whether $n$ is divisible by $4$ or not.
For $Z=\IZ_{2}\times\IZ_{2}$, we get $\<\cow{1},\cow{n}\>=\half{1}$
so that $\ellb=2$. To determine $\ellf$, notice that if there
exists a commutator map $\omega$ on $\lambdaz=\weight$ satisfying
\eqref{commutator} with $\ell=1$, the fundamental level
of $\Spin_{2n}/\IZ_{2}^{\pm}$ is one and therefore $n$
is divisible by 4. Conversely, if $4|n$, the form
$\wt\omega(\alpha\oplus p\cow{1}\oplus q\cow{n},
           \beta\oplus p'\cow{1}\oplus q'\cow{n})=
 i^{pq'-qp'}(-1)^{\<\alpha,\beta\>+\<\alpha,p'\cow{1}+
 q'\cow{n}\>+\<p\cow{1}+q\cow{n},\beta\>}$
defined on $\root\oplus\cow{1}\IZ\oplus\cow{n}\IZ$ descends to
a suitable form on
$\weight=\root\oplus\cow{1}\IZ\oplus \cow{n}\IZ/
 (-2\cow{1}\oplus 2\cow{1}\oplus 0)\IZ+(-2\cow{n}\oplus 0\oplus 2\cow{n})\IZ$.

${\bf E_{6}}$\\
$E_{6}$ is simply--laced and $\weight/\root\cong\IZ_{3}$ is generated
by any of its non--zero elements and therefore by
$\cow{6}=\theta_{5}-\frac{1}{3}(\theta_{6}+\theta_{7}-\theta_{8})$
corresponding to the special root $\alpha_{6}=-\theta_{4}+\theta_{5}$.
Since $\<\cow{6},\cow{6}\>=\frac{4}{3}$, $\ell_{b}=3$ and $\ellf=1$.

${\bf E_{7}}$\\
$E_{7}$ is simply--laced and $\weight/\root\cong\IZ_{2}$ is generated
by $\cow{7}=\theta_{6}-\half{1}(\theta_{7}-\theta_{8})$ corresponding
to the unique special root $\alpha_{7}=-\theta_{5}+\theta_{6}$. Since
$\<\cow{7},\cow{7}\>=\half{3}$, $\ell_{b}=\ellf=2$ \halmos

\ssection{The action of $\lzg$ on the positive energy dual of $LG$}
\label{se:disc loops}
%==================================================================

We show in this section that the category $\pl$ of positive energy
representations of $LG$ at a given level $\ell$ is closed under
conjugation by elements of $\lzg$. We also identify the corresponding
abstract action of $Z\cong\lzg/LG$ on the alcove of $G$ parametrising
the irreducibles in $\pl$ with the geometric one obtained by realising
$Z$ as a distinguished subgroup of the automorphisms of the extended
Dynkin diagram of $G$. We begin by studying the latter.

\ssubsection{Geometric action of $\mathbf{Z(G)}$ on the level $\mathbf{\ell}$ alcove}
\label{ss:Z on alcove}
%------------------------------------------------------------------------------------

This subsection is essentially an expanded version of \cite[ch. VI, \S 2.3]{Bou}.
The notation follows that of section \ref{se:lattices}. Denote $-\theta$
by $\alpha_{0}$, then

\begin{lemma}
For any special root $\alpha_{i}$, the set
$\Delta_{i}=\Delta\backslash\{\alpha_{i}\}\cup\{\alpha_{0}\}$ is a basis
of $R$ with highest root $-\alpha_{i}$ and dual basis
\begin{align}
\cowp{0}&=-\cow{i} \label{dual one}\\
\cowp{j}&= \cow{j}-\<\theta,\cow{j}\>\cow{i} \label{dual two}
\end{align}
\end{lemma}
\proof
Let $x\in\tc$, then, by \eqref{highest root}
\begin{equation}
 x=\sum\<x,\cow{j}\>\alpha_{j}
  =\<x,\cow{i}\>\theta+\sum_{j\neq i}
   (\<x,\cow{j}\>-\<x,\cow{i}\>\<\theta,\cow{j}\>)\alpha_{j}
\end{equation}
so that $\Delta_{i}$ is a vector space basis of $\tc$ with dual basis given
by \eqref{dual one}--\eqref{dual two}.
If $\beta$ is a positive root, then either $\<\beta,\cow{i}\>=0$, in which
case $\<\beta,\cowp{0}\>$ and $\<\beta,\cowp{j}\>$ are all non--negative,
or $\<\beta,\cow{i}\>=1$ since $\cow{i}$ is a minimal dominant coweight.
In the latter case $\<\beta,\cowp{0}\>=-1$ and
$\<\beta,\cowp{j}\>=\<\beta-\theta,\cow{j}\>\leq 0$.
Thus, $\Delta_{i}$ is a basis of $R$.
Next, for any $\beta\in R$,
$\<-\alpha_{i}-\beta,\cowp{0}\>=1+\<\beta,\cow{i}\>\geq 0$ since $\cow{i}$
is minimal. Moreover, for $j\neq i$
\begin{equation}\label{positive}
 \<-\alpha_{i}-\beta,\cowp{j}\>=
 \<\theta,\cow{j}\>-\<\beta,\cow{j}\>+\<\theta,\cow{j}\>\<\beta,\cow{i}\>=
 \<\theta-\beta,\cow{j}\>+\<\theta,\cow{j}\>\<\beta,\cow{i}\>
\end{equation}
The above is clearly non--negative if $\<\beta,\cow{i}\>\geq 0$. If, on
the other hand, $\<\beta,\lambda_{i}\vvee\>=-1$ then $\beta$ is negative
and \eqref{positive} is equal to $-\<\beta,\lambda_{j}\vvee\>\geq 0$. Thus,
$-\alpha_{i}$ is the highest root relative to $\Delta_{i}$ \halmos

\begin{proposition}\label{centre weyl}
Let $\deltabar=\Delta\cup\{\alpha_{0}\}$. Then, for any special root
$\alpha_{i}$, there exists a unique
\begin{equation}
w_{i}\in W_{0}=\{w\in W|\thinspace w\deltabar=\deltabar\}
\end{equation}
such that $w_{i}\alpha_{0}=\alpha_{i}$. The resulting map
$\imath:Z(G)\rightarrow W_{0}$ obtained by identifying $Z(G)\setminus\{1\}$
with the set of special roots is a group isomorphism.
\end{proposition}
\proof
The existence of $w_{i}$ follows from the previous lemma since $W$ acts
transitively on the set of basis of $R$ and maps highest roots to highest
roots. $w_{i}$ is unique because an element $w\in W_{0}$ is determined by
$w\alpha_{0}$. Indeed, if $w_{1}\alpha_{0}=\alpha_{j}=w_{2}\alpha_{0}$,
then $w_{2}^{-1}w_{1}$ is a permutation of $\Delta$ and is therefore the
identity since $W$ acts simply on basis. $\imath$ is injective because
$w_{i}\alpha_{0}=\alpha_{i}$. Let now $w\in W_{0}$. We claim that
$\alpha_{i}=w\alpha_{0}$ is a special root. It then follows by uniqueness
that $w=w_{i}$ and therefore that $\imath$ is surjective. To see this, we
apply $w$ to \eqref{highest root} and get
$-\alpha_{i}=\sum_{j}m_{j}w\alpha_{j}$ while at the same time
$-\alpha_{i}=m_{i}^{-1}(\alpha_{0}+\sum_{j\neq i}m_{j}\alpha_{j})$.
Equating the coefficients of $\alpha_{0}$, we get $m_{i}=1$. To prove
that $\imath$ is a homomorphism, let $\alpha_{i}$ and $\alpha_{j}$ be special
roots. Then either $w_{i}w_{j}=1$ or $w_{i}w_{j}=w_{k}$ where $\alpha_{k}$
is another special root. In the former case, $w_{i}\alpha_{j}=\alpha_{0}$
and therefore, by \eqref{dual one}
\begin{equation}\label{ccompo}
w_{i}\cow{j}=\cowp{0}=-\cow{i}
\end{equation}
so that $\cow{j}=-\cow{i}$ mod $\coroot$ since $W$ leaves $\coroot$--cosets
invariant. In the latter, $w_{i}\alpha_{j}=\alpha_{k}$ and therefore, using
\eqref{dual two}
\begin{equation}\label{compo}
w_{i}\cow{j}=\cowp{k}=\cow{k}-\<\theta,\cow{k}\>\cow{i}=\cow{k}-\cow{i}
\end{equation}
whence $\cow{i}+\cow{j}=\cow{k}$ mod $\coroot$ \halmos\\

The following is well--known and often rediscovered \cite{OT,Ga}

\begin{corollary}
$Z(G)$ is canonically isomorphic to the group of automorphisms of the
extended Dynkin diagram of $G$ induced by Weyl group elements.
\end{corollary}

For any $\ell\in\IN$, recall that the level $\ell$ alcove is the set
defined by
\begin{equation}
\al=\{\lambda\in\weight|
    \<\lambda,\alpha_{i}\>\geq 0,\medspace\<\lambda,\theta\>\leq\ell\}
\end{equation}

\begin{proposition}\label{geometric action}
For any $\ell\in\IN$, there is a canonical action of $Z(G)$ on the level
$\ell$ alcove $\al$ given by
\begin{equation}\label{eq:action of Z}
z\longrightarrow A_{i}=\tau(\ell\cow{i})w_{i}
\end{equation}
where $i$ is the label of the special root corresponding to $z$ via lemma
\ref{special roots}, $\tau$ denotes translation and $w_{i}=\imath(z)$ corresponds
to $z$ via proposition \ref{centre weyl}.
\end{proposition}
\proof
If $\lambda\in\al$ then for $j\neq i$,
$\<A_{i}\lambda,\alpha_{j}\>=\<\lambda,w_{i}^{-1}\alpha_{j}\>\geq 0$ since
$w_{i}^{-1}\alpha_{j}\neq \alpha_{0}$. On the other hand,
$\<A_{i}\lambda,\alpha_{i}\>=\ell+\<\lambda,\alpha_{0}\>\geq 0$. Finally,
$\<A_{i}\lambda,\theta\>=\ell-\<\lambda,w_{i}^{-1}\theta\>\leq\ell$ so that
the $A_{i}$ leave $\al$ invariant. Next,
$A_{i}A_{j}=\tau(\ell(\cow{i}+w_{i}\cow{j}))w_{i}w_{j}$. If $w_{i}w_{j}=1$,
we get by \eqref{ccompo} and the previous proposition $A_{i}A_{j}=1$. If, on
the other hand $w_{i}w_{j}=w_{k}$, \eqref{compo} yields $A_{i}A_{j}=A_{k}$
\halmos\\

\remark The explicit action of $Z(G)$ on $\al$ for all classical groups is
given in \S \ref{ss:Z on alcove 2}.

\ssubsection{Positive energy representations of $\mathbf{LG}$}
\label{ss:per of LG}
%-------------------------------------------------------------

We outline the classification of positive energy representations
of $LG$ following \cite{Wa} to which we refer for more details.
Let $\pi$ be a projective unitary representation of $LG\rtimes\rot$
on a complex Hilbert space $\H$, \ie a strongly continuous homomorphism
\begin{equation}
\pi:LG\rtimes\rot\longrightarrow PU(\H)=U(\H)/\T
\end{equation}
Over $\rot$, $\pi$ lifts to a unitary representation which we denote by
the same symbol. By definition, $\H$ is of positive energy if
\begin{equation}
\H=\bigoplus_{n\geq n_{0}}\H(n)
\end{equation}
where each $\H(n)=\{\xi\in\H|\pi(R_{\theta})\xi=e^{in\theta}\xi\}$ is
finite--dimensional. The lift is unique up to multiplication by a character
of $\rot$ and we normalise it by choosing $n_{0}=0$ and $\H(0)\neq\{0\}$.\\
%Since $G$ is simple, connected and simply--connected, the restriction of $\pi$
%to $G$ lifts uniquely to a unitary representation which commutes with $\rot$.
%In particular, each $\H(n)$ is a finite--dimensional $G$--module.\\

The classification of positive energy representation is obtained via
the associated infinitesimal action of the Lie algebra of $\g$--valued
trigonometric polynomials $\lpol\g\subset L\g$ in the following way. Consider
the subspace $\hfin\subset\H$ of finite energy vectors for $\rot$, that
is the algebraic direct sum of the $\H(n)$. The latter is a core for the
normalised self--adjoint generator of rotations which we denote by $d$.
Thus
\begin{xalignat}{3}
\left.d\right|_{\H(n)}&=n&&\text{and}&\pi(R_{\theta})&=e^{i\theta d}
\end{xalignat}
For any $X\in\lpol\g$, the one--parameter projective group $\pi(\exp_{LG}(tX))$
possesses a continuous lift to $U(\H)$, unique up to multiplication by a
character of $\IR$. It is therefore given, via Stone's theorem by
$e^{t\pi(X)}$ where $\pi(X)$ is a skew--adjoint operator determined up to
an additive constant.

\begin{theorem}[Wassermann]\label{th:core}
The subspace $\Hfin$ of finite energy vectors is an invariant core for the
operators $\pi(X)$, $X\in\lpol\g$. The operators $\pi(X)$ may be chosen
uniquely so as to satisfy $[d,\pi(X)]=i\pi(\dot X)$ on $\Hfin$ and then
$X\rightarrow\pi(X)$ gives a projective representation of $\lpol\g$ on
$\Hfin$ such that
\begin{equation}\label{eq:KM reln}
[\pi(X),\pi(Y)]=\pi([X,Y])+i\ell B(X,Y)
\end{equation}
where $B(X,Y)=\int_{0}^{2\pi}\<X,\dot Y\>\frac{d\theta}{2\pi}$ and $\ell$
is a non--negative integer called the level of $\H$.
\end{theorem}

We denote the restriction of the operators $\pi(X)$, $X\in\lpol\g$ to $\Hfin$
by the same symbol and extend $\pi$ to a projective representation of
$\lpol\gc$ on $\hfin$ satisfying \eqref{eq:KM reln} as well as the formal
adjunction property $\pi(X)^{*}=-\pi(\overline{X})$.
The operators $\pi(X)$ and $d$ then give rise to a unitarisable
representation of the Kac--Moody algebra $\wh\gc$ at level $\ell$ such that
$d$ is diagonal with finite--dimensional eigenspaces and spectrum in $\IN$.
Such representations split into a direct sum of irreducibles, each of which
is an integrable highest weight representation, that is a module
generated over the enveloping algebra $\U(\g_{\leq})$ by a vector $v$
uniquely determined by the requirement that it be 
annihilated by $\g_{\geq}$ and diagonalises the action of $T\rtimes\rot$.
Here $\g_{\leq}$ (resp. $\g_{\geq}$) is the nilpotent Lie algebra spanned
by the $x(n)$ with $n<0$ (resp. $n>0$) and $x\in\gc$ or $n=0$ and $x$ lying
in a negative (resp. positive) root space of $\gc$. Thus, for any $h\in\tc$
\begin{align}
d v&=n v\\
\pi(h) v&=\lambda(h) v
\end{align}
for some $n\in\IN$ and dominant integral weight $\lambda$ of $G$ which
satisfies $\<\lambda,\theta\>\leq\ell$. The pair $(\ell,\lambda)$ classifies
the integrable representation uniquely \cite{Ka}. The finite collection $\al$
of dominant integral weights of $G$
satisfying $\<\lambda,\theta\>\leq\ell$ is called the level $\ell$ alcove
of $G$.\\

The classification of a positive energy representation $(\pi,\H)$ as
an $LG$--module is equivalent as that of $\hfin$ as an $\lpol\g$--module.
In particular, $\H$ is topologically irreducible under $LG$ iff $\hfin$
is irreducible under $\lpol\g$ and is then uniquely determined by its
level $\ell$ and highest weight $\lambda$. Moreover, for any $\ell\in\IN$
and $\lambda\in\al$, there exists a unique unitarisable highest weight
$\wh{\gc}$--module with level $\ell$ and highest weight $\lambda$ \cite{Ka}.
The corresponding action of $\lpol\g$ may then be exponentiated to yield
a positive energy representation of $LG$ on its Hilbert space completion
\cite{GoWa,TL2}.

\ssubsection{Action of $\bf{\lzg}$ on positive energy representations of $LG$}
\label{ss:Z on per}
%-----------------------------------------------------------------------------

The group of discontinuous loops $\lzg$ acts on $LG$ by conjugation and
therefore on the irreducible projective unitary representations of $LG$
by $\zeta_{*}\pi(\gamma)=\pi(\zeta^{-1}\gamma\zeta)$, $\zeta\in\lzg$.
The following result shows that positive energy ones are stable under
this action

\begin{proposition}\label{pr:stability}
If $(\pi,\H)$ is an irreducible positive energy representation of
$LG$ and $\zeta\in\lzg$, the conjugated representation $\zeta_{*}\pi$
is of positive energy.
\end{proposition}
\proof
Any intertwining action of $\rot$ for $\zeta_{*}\pi$, whether of
positive energy or not is necessarily given by 
\begin{equation}\label{twisted rotations}
\theta\longrightarrow V_{\theta}=
\pi(\zeta^{-1}\zeta_{\theta})U_{\theta}
\end{equation}
where $\zeta_{\theta}(x)=\zeta(x-\theta)$ and $U_{\theta}$ is the
positive energy intertwining action of $\rot$ for $\pi$. Indeed,
$\zeta^{-1}\zeta_{\theta}\in LG$ and $V_{\theta}$ yields a projective
action of $\rot$ satisfying
$V_{\theta}\zeta_{*}\pi(\gamma)V_{\theta}^{*}=\zeta_{*}\pi(\gamma_{\theta})$.
Moreover, if $V_{\theta}^{i}$, $i=1,2$ are two actions of $\rot$
intertwining $\zeta_{*}\pi$, then
$W_{\theta}=(V_{\theta}^{1})^{*}V_{\theta}^{2}$ commutes
projectively with $LG$ and the following holds in $U(\H)$
\begin{equation}
W_{\theta}\pi(\gamma)W_{\theta}^{*}\pi(\gamma)^{*}=\chi(\gamma,\theta)
\end{equation}
where $\chi(\gamma,\theta)\in\T$ depends multiplicatively on either
variable. Since $LG$ is perfect, \cite[prop. 3.4.1.]{PS}, $\chi\equiv 1$
and by Shur's lemma $W_{\theta}=1$ in $PU(\H)$.\\

Thus, $\zeta_{*}\pi$ is of positive energy iff $V_{\theta}$ is a
positive energy representation of $\rot$. It is sufficient to check
this for a given set of representatives of $LG$--cosets in $\lzg$
and therefore for the discontinuous loops
$\zeta_{\mu}(\theta)=\exp_{T}(-i\theta\mu)$, $\mu\in\lambdaz$. If
$\mu\in\coroot$, then $\zeta_{\mu}\in LG$ and the action of $\rot$
given by \eqref{twisted rotations} may be rewritten as
$\pi(\zeta^{-1}_{\mu}{\zeta_{\mu}}_{\theta})U_{\theta}=
 \pi(\zeta_{\mu})^{*}U_{\theta}\pi(\zeta_{\mu})$ which is of positive
energy. The general case $\mu\in\lambdaz$ is settled by the
following simple observation. Notice that
$\zeta_{\mu}^{-1}{\zeta_{\mu}}_{\theta}=
 \zeta_{\mu}(-\theta)=\exp_{T}(i\theta\mu)$
since $\zeta_{\mu}$ is a homomorphism and write $\mu$ as a convex
combination of elements in the coroot lattice,
$\mu=\sum_{i=1}^{m} t_{i}\mu_{i}$, $t_{i}\in(0,1]$, $\sum_{i} t_{i}=1$,
$\mu_{i}\in\coroot$. Since $\pi$ lifts to a unitary representation
$\wt\pi$ over $G\times\rot\supset T\times\rot$, 
\begin{equation}\label{convex}
\wt\pi(\exp_{T}(i\theta\mu))\wt\pi(R_{\theta})=
\prod_{j}\wt\pi(\exp_{T}(i\theta t_{j}\mu_{j}))\wt\pi(R_{t_{j}\theta})
\end{equation}
is a lift of $\pi(\zeta_{\mu}^{-1}{\zeta_{\mu}}_{\theta})U_{\theta}$
and the product of $m$ commuting representations of $\IR$ which
by our previous argument are of positive energy. It follows that
$\pi(\zeta_{\mu}^{-1}{\zeta_{\mu}}_{\theta})U_{\theta}$ is of
positive energy \halmos

\begin{proposition}\label{conjugated action}
Let $(\pi,\H)$ be a positive energy representation of $LG$ of level
$\ell$ and $\zeta\in\lzg$. Then,
\begin{enumerate}
\item $\zeta_{*}\pi$ is of level $\ell$.
\item If $\zeta(\phi)=\exp_{T}(-i\phi\mu)$ is the discontinuous loop
corresponding to $\mu\in\lambdaz$, the subspaces of finite energy
vectors of $\pi$ and $\zeta_{*}\pi$ coincide.
\item If $\Ad(\zeta)\lpol\g=\lpol\g$ and the finite energy subspaces of
$\pi$ and $\zeta_{*}\pi$ coincide, the conjugated action of $\lpol\g$
on $\hfin$ is given by
\begin{equation}\label{twisted lg}
\zeta_{*}\pi(X)=
\pi(\zeta^{-1} X\zeta)+
i\ell\int_{0}^{2\pi}\<\dot\zeta\zeta^{-1},X\>\frac{d\theta}{2\pi}
\end{equation}
\end{enumerate}
\end{proposition}
\proof
It is sufficient to check (i) for the discontinuous loops $\zeta_{\mu}$,
$\mu\in\lambdaz$. This will be done in the course of the proof of (iii).

(ii)
It was remarked in the proof of the previous proposition that the conjugated
action of rotations \eqref{twisted rotations} corresponding to $\zeta_{\mu}$
is given by
\begin{equation}
 \pi(\zeta^{-1}_{\mu}{\zeta_{\mu}}_{\theta})U_{\theta}=
 \pi(\exp_{T}(i\theta\mu))U_{\theta}
\end{equation}
which commutes with the original action of $\rot$ given by $U_{\theta}$.
Since both are of positive energy, their finite energy subspaces coincide.

(iii)
Let $h\in\IN$ be the level of $\zeta_{*}\pi$ and denote by $\pi$ and
$\zeta_{*}\pi$ the projective representations of $\lpol\g$ on $\hfin$
given by theorem \ref{th:core}, so that
\begin{align}
[\pi(X),\pi(Y)]&=\pi([X,Y])+i\ell B(X,Y)\\
[\zeta_{*}\pi(X),\zeta_{*}\pi(Y)]&=\zeta_{*}\pi([X,Y])+ihB(X,Y)
\label{new level}
\end{align}
Evidently,
$\zeta_{*}\pi(X)=\pi(\zeta^{-1}X\zeta)+iF(X)$ for some $F(X)\in\IR$ since
$\zeta_{*}\pi(\exp_{LG}(X))= \pi(\exp_{LG}(\zeta^{-1}X\zeta))=
 e^{\pi(\zeta^{-1}X\zeta)}$ in $PU(\H)$. It follows, by \eqref{eq:zeta on B}
that
\begin{equation}
\begin{split}
[\zeta_{*}\pi(X),\zeta_{*}\pi(Y)]
&=
[\pi(\zeta^{-1}X\zeta),\pi(\zeta^{-1}Y\zeta)]\\
&=
\pi(\zeta^{-1}[X,Y]\zeta)+i\ell B(\zeta^{-1}X\zeta,\zeta^{-1}Y\zeta)\\
&=
\zeta_{*}\pi([X,Y])-iF([X,Y])+i\ell B(X,Y)+i\ell\int_{0}^{2\pi}
\<\dot\zeta\zeta^{-1},[X,Y]\>
\frac{d\theta}{2\pi}
\end{split}
\end{equation}
Since $hB$ and $\ell B$ lie in the same cohomology class iff $h=\ell$,
we find by equating the above with \eqref{new level} that the level of
$\zeta_{*}\pi$ is $\ell$. Moreover, \eqref{twisted lg} holds since
$[\lpol\g,\lpol\g]=\lpol\g$ \halmos

\begin{theorem}\label{th:Z on per}
Let $(\pi,\H)$ be an irreducible positive energy representation of $LG$
of level $\ell$ and highest weight $\lambda$ and $\zeta\in\lzg$. Then,
the conjugated representation
$\zeta_{*}\pi(\gamma)=\pi(\zeta^{-1}\gamma\zeta)$ on $\H$ is of positive
energy, level $\ell$ and highest weight $\zeta\lambda$ where the notation
refers to the geometric action of $Z(G)\cong L_{Z(G)}G/LG$ on the level
$\ell$ alcove defined by proposition \ref{geometric action}.
\end{theorem}
\proof
It suffices to prove the result for a given choice of representatives
of $LG$--cosets in $\lzg$. Let $z\in Z\backslash\{1\}$ correspond
to the special root $\alpha_{j}$ by lemma \ref{special roots} and consider
the discontinuous loop $\zeta=\zeta_{\cow{j}}w_{j}$ where $\cow{j}$ is
the associated fundamental coweight and $w_{j}\in G$ a representative
of the Weyl group element corresponding to $z$ by proposition
\ref{centre weyl}. Since $G$ commutes $\rot$, the subspace of finite
energy vectors of $\pi$ coincides with that of ${w_{j}}_{*}\pi$ and,
by the previous proposition, with that of $\zeta_{*}\pi$. We may
therefore compare the infinitesimal actions of $\lpol\gc$ corresponding
to $\pi$ and $\zeta_{*}\pi$ on $\Hfin$.\\
 
If $\alpha$ is a root and $x_{\alpha}\in\g_{\alpha}$, then
$[\cow{j},x_{\alpha}]=\<\cow{j},\alpha\>x_{\alpha}$. Since
$\zeta_{\cow{j}}(\theta)=\exp_{T}(-i\cow{j}\theta)$, this gives
\begin{equation}
 \zeta_{\cow{j}}^{-1}x_{\alpha}(n)\zeta_{\cow{j}}(\theta)=
 x_{\alpha}\otimes e^{i\theta(n+\<\cow{j},\alpha\>)}=
 x_{\alpha}(n+\<\cow{j},\alpha\>)(\theta)
\end{equation}
Therefore, up to a non--zero multiplicative constant
\begin{equation}\label{new e alpha}
\zeta_{*}\pi(e_{\alpha}(n))=
\pi(e_{w_{j}^{-1}\alpha}(n+\<\cow{j},\alpha\>))
\end{equation}
since no additional term arises from \eqref{twisted lg} because
$\dot\zeta\zeta^{-1}=-i\cow{j}$ lies in $\tc$ which is orthogonal to
$\g_{\alpha}$. If, on the other hand $h\in\tc$, then
$\zeta^{-1}h(n)\zeta=w_{j}^{-1}h(n)$ and \eqref{twisted lg} reads
\begin{equation}\label{new hn}
\zeta_{*}\pi(h(n))=\pi(w_{j}^{-1}h(n))+\ell\delta_{n,0}\<h,\cow{j}\>
\end{equation}
Let $\Omega\in\hfin$ be the highest weight vector for $\zeta_{*}\pi$.
We claim that, up to a scalar factor, $\Omega=\Upsilon$, the highest
weight vector for $\pi$. To see this, recall that $\Omega$ is the unique
element of $\hfin$ annihilated by the subalgebra $\g_{\geq}$ spanned
by the $x(n)$, $x\in\gc$ and $n>0$ and the $x_{\alpha}(0)$ with
$\alpha>0$. $\g_{\geq}$ is generated by the elements corresponding
the simple affine roots, namely $e_{\alpha_{i}}(0)$ and
$e_{\alpha_{0}}(1)$ where $\alpha_{0}=-\theta$.
Recalling from proposition \ref{centre weyl} that $w_{j}^{-1}$ acts
as a permutation of $\deltabar=\{\alpha_{0},\ldots,\alpha_{n}\}$ and
maps $\alpha_{j}$ to $\alpha_{0}$, we get, using \eqref{new e alpha}
\begin{align}
\zeta_{*}\pi(e_{\alpha_{k}}(0))
&=
\left\{\begin{array}{ll}
\pi(e_{w_{j}^{-1}\alpha_{k}}(0))&\text{if $k\neq j$}\\
\pi(e_{\alpha_{0}}(1))          &\text{if $k=j$}
\end{array}\right.\\
\zeta_{*}\pi(e_{\alpha_{0}}(1))
&=\pi(e_{w_{j}^{-1}\alpha_{0}}(0))
\end{align}
whence $\Omega=\Upsilon$. To find the weight of $\Omega$ and therefore
the highest weight of $\zeta_{*}\pi$, we use \eqref{new hn} and the fact
that $\pi(h(0))\Upsilon=\<\lambda,h\>\Upsilon$ whenever $h\in\tc$ so that
\begin{equation}
\zeta_{*}\pi(h(0))\Omega=
\<h,w_{j}\lambda+\ell\cow{j}\>\Omega=
\<h,\zeta\lambda\>\Omega
\end{equation}
\halmos

The proof of the above theorem has the following useful
\begin{corollary}\label{automorphism}
Let $z\in Z(G)\backslash\{1\}$ and $\cow{j}\in\coweight$ and $w_{j}\in G$ the
fundamental coweight and representative of the Weyl group element corresponding
to $z$ by lemma \ref{special roots} and proposition \ref{centre weyl}
respectively. Then, conjugation by $\zeta=\zeta_{\cow{j}}w_{j}$ induces an
automorphism of $\lpol\gc$ preserving its triangular decomposition.
\end{corollary}

\ssubsection
{Appendix : explicit action of $\bf{Z(G)}$ on the level $\mathbf{\ell}$ alcove}
\label{ss:Z on alcove 2}
%------------------------------------------------------------------------------

We describe explicitly the action of $Z(G)$ given by proposition
\ref{geometric action} for all classical groups, using the tables
\cite[Tables I--IV]{Bou}. For each special root $\alpha_{i}$, we
denote the corresponding element of $W_{0}\subset W$ by $w_{i}$
and note that the fundamental coweight $\cow{i}$ and weight
$\lambda_{i}$ coincide since $\alpha_{i}$ is long.
We let moreover $\theta_{i}$, $i=1\ldots n$ be the standard
basis in $\IR^{n}$ and $I$ the lattice $\bigoplus_{i}\theta_{i}
\cdot\IZ$. Unless otherwise stated, the basic inner product is
the standard one on $\IR^{n}$.\\

{\bf SU$_{\bf n}$, ${\bf n\geq 2}$}\\
%-----------------------------------------------------
simple roots : $\alpha_{i}=\theta_{i}-\theta_{i+1}$, $i=1\ldots n-1$.\\
highest root : $\theta=\theta_{1}-\theta_{n}=\alpha_{1}+\cdots+\alpha_{n-1}$.\\
minimal dominant coweights :
$\cow{i}=\theta_{1}+\cdots+\theta_{i}-\frac{i}{n}\sum_{j}\theta_{j}$,
$i=1\ldots n-1$.\\
Weyl group : $\mathfrak S_{n}$ acting by permutation of the $\theta_{i}$.\\
$W_{0}$ : $w_{k}$ is the cyclic permutation
$(\theta_{1}\cdots\theta_{n})^{k}=(\alpha_{0}\cdots\alpha_{n-1})^{k}$.\\
level $\ell$ alcove :
$\al=\{\mu\in I|\thinspace\mu_{1}\geq\cdots\geq\mu_{n},\thickspace
                              \mu_{1}-\mu_{n}\leq\ell\}/(\sum_{j}\theta_{j})$.\\
action of the centre :
$A_{k}(\mu_{1},\ldots,\mu_{n})=
 (\ell+\mu_{n+1-k},\ldots,\ell+\mu_{n},\mu_{1},\ldots,\mu_{n-k})$.\\

{\bf Spin$_{\bf 2n+1}$, ${\bf n\geq 2}$}\\
%---------------------------------------
simple roots :
$\alpha_{i}=\theta_{i}-\theta_{i+1}$, $i=1\ldots n-1$ and $\alpha_{n}=\theta_{n}$.\\
highest root :
$\theta=\theta_{1}+\theta_{2}=\alpha_{1}+2(\alpha_{2}+\cdots+\alpha_{n})$.\\
minimal dominant coweight : $\cow{1}=\theta_{1}$.\\
Weyl group :
$\mathfrak S_{n}\ltimes\IZ_{2}^{n}$ acts by permutations and sign changes
of the $\theta_{i}$.\\
$W_{0}$ :
$w_{1}$ is the sign change $\theta_{1}\rightarrow-\theta_{1}$ permuting
$\alpha_{0}$ and $\alpha_{1}$.\\
level $\ell$ alcove :
$\al=\{\mu\in I+\half{1}(\theta_{1}+\ldots+\theta_{n})\IZ|
     \thinspace\mu_{1}\geq\cdots\geq\mu_{n}\geq 0,\medspace \mu_{1}+\mu_{2}\leq\ell\}$.\\
action of the centre :
$A_{1}(\mu_{1},\mu_{2},\ldots,\mu_{n})=(\ell-\mu_{1},\mu_{2},\ldots,\mu_{n})$.\\

{\bf Sp$_{\bf n}$, ${\bf n\geq 2}$}\\
%----------------------------------
simple roots :
$\alpha_{i}=\theta_{i}-\theta_{i+1}$, $i=1\ldots n-1$ and $\alpha_{n}=2\theta_{n}$.\\
highest root : $\theta=2\theta_{1}=2(\alpha_{1}+\cdots+\alpha_{n-1})+\alpha_{n}$.\\
basic inner product : half the standard one on $\IR^{n}$.\\
minimal dominant coweight : $\cow{n}=\theta_{1}+\cdots+\theta_{n}$.\\
Weyl group : $\mathfrak S_{n}\ltimes\IZ_{2}^{n}$ acts by permutations
and sign changes of the $\theta_{i}$.\\
$W_{0}$ : $w_{n}$ is the transformation $\theta_{i}\rightarrow -\theta_{n+1-i}$.\\
level $\ell$ alcove :
$\al=\{\mu\in I|\thinspace \ell\geq\mu_{1}\geq\cdots\geq\mu_{n}\geq 0\}$.\\
action of the centre :
$A_{n}(\mu_{1},\ldots,\mu_{n})=(\ell-\mu_{n},\ldots,\ell-\mu_{1})$.\\

{\bf Spin$_{\bf 2n}$, ${\bf n\geq 3}$}\\
%-----------------------------------------------------
simple roots : $\alpha_{i}=\theta_{i}-\theta_{i+1}$, $i=1\ldots n-1$ and
$\alpha_{n}=\theta_{n-1}+\theta_{n}$.\\
highest root :
$\theta=
\theta_{1}+\theta_{2}=
\alpha_{1}+2(\alpha_{2}+\cdots+\alpha_{n-2})+\alpha_{n-1}+\alpha_{n}$.\\
minimal dominant coweights :
$\lambda_{1}  =\theta_{1}$,
$\lambda_{n-1}=\half{1}(\theta_{1}+\cdots+\theta_{n-1}-\theta_{n})$,
$\lambda_{n}  =\half{1}(\theta_{1}+\cdots+\theta_{n})$.\\
Weyl group : $\mathfrak S_{n}\ltimes\IZ_{2}^{n-1}$ acts by permutations
and even numbers of sign changes of the $\theta_{i}$.\\
$W_{0}$ : $w_{1}$ is the sign change $\theta_{1}\rightarrow -\theta_{1}$,
$\theta_{n}\rightarrow -\theta_{n}$ and permutes $\{\alpha_{0},\alpha_{1}\}$
and $\{\alpha_{n-1},\alpha_{n}\}$.
For $n$ even, $w_{n-1}$ is given by $\theta_{i}\rightarrow -\theta_{n+1-i}$,
$2\leq i\leq n-1$ and $\theta_{1}\leftrightarrow\theta_{n}$ and permutes
$\{\alpha_{0},\alpha_{n-1}\}$ and $\{\alpha_{1},\alpha_{n}\}$ while
$w_{n}$ is given by $\theta_{i}\rightarrow -\theta_{n+1-i}$ and permutes
$\{\alpha_{k},\alpha_{n-k}\}$.
For $n$ odd, $w_{n-1}$ is given by $\theta_{1}\rightarrow\theta_{n}$ and
$\theta_{i}\rightarrow -\theta_{n+1-i}$, $i=2\ldots n$ and acts as the cyclic
permutation
$\begin{pmatrix}\alpha_{1}&\alpha_{n}&\alpha_{0}&\alpha_{n-1}\end{pmatrix}$
while $w_{n}$ is given by $\theta_{i}\rightarrow -\theta_{n+1-i}$, $i=1\ldots n-1$
and $\theta_{n}\rightarrow\theta_{1}$ and acts as
$\begin{pmatrix}\alpha_{1}&\alpha_{n}&\alpha_{0}&\alpha_{n-1}\end{pmatrix}^{-1}$.\\
level $\ell$ alcove :
$\al=\{\mu\in I+\half{1}(\theta_{1}+\ldots+\theta_{n})\IZ|
     \thinspace\mu_{1}\geq\cdots\geq\mu_{n-1}\geq|\mu_{n}|,
     \medspace\mu_{1}+\mu_{2}\leq\ell\}$.\\
action of the centre :
\begin{align}
A_{1}(\mu_{1},\ldots,\mu_{n})&=(\ell-\mu_{1},\mu_{2},\ldots,\mu_{n-1},-\mu_{n}).\\
A_{n-1}(\mu_{1},\ldots,\mu_{n})&=
\begin{cases}
(\half{\ell}+\mu_{n},\half{\ell}-\mu_{n-1},\ldots,\half{\ell}-\mu_{2},-\half{\ell}+\mu_{1})
&\text{$n$ even}\\[1.2 ex]
(\half{\ell}-\mu_{n},\ldots,\half{\ell}-\mu_{2},-\half{\ell}+\mu_{1})
&\text{$n$ odd}
\end{cases}\label{eq:n-1}\\
A_{n}(\mu_{1},\ldots,\mu_{n})&=
\begin{cases}
(\half{\ell}-\mu_{n},\ldots,\half{\ell}-\mu_{1})
&\text{$n$ even}\\[1.2 ex]
(\half{\ell}+\mu_{n},\half{\ell}-\mu_{n-1},\ldots,\half{\ell}-\mu_{1})
&\text{$n$ odd}\label{eq:n}
\end{cases}
\end{align}

\ssection{The Mackey obstruction corresponding to $LG\subset\lzg$}
\label{se:Mackey obstruction}
%=================================================================

We determine below the central extensions of $\lzg$ corresponding to
positive energy representations which remain irreducible when restricted
to $LG$. More precisely, let $\H$ be an irreducible, level $\ell$ positive
energy representation of $LG$ the isomorphism class of which is invariant
under $\lzg$. $\H$ gives rise to a projective action $\pi$ of $\lzg$
extending that of $LG$ and therefore by pull--back of
\begin{equation}
1\rightarrow\T\rightarrow U(\H)\rightarrow PU(\H)\rightarrow 1
\end{equation}
to a central extension $\wt\lzg$ of $\lzg$. Since the restriction
of $\wt\lzg$ to $LG$ is smooth and of level $\ell$
\cite[chap. II, \S 2.4]{TL}, $\wt\lzg$ is smooth and therefore falls
within the classification of section \ref{se:extensions of lzg}.
In particular, it is uniquely determined by its commutator map
\begin{equation}\label{eq:yet another}
\omega(\lambda,\mu)=
\pi(\zeta_{\lambda})\pi(\zeta_{\mu})
\pi(\zeta_{\lambda})^{*}\pi(\zeta_{\mu})^{*}
\end{equation}
which satisfies
\begin{equation}\label{triviality}
\omega(\alpha,\mu)=(-1)^{\ell\<\alpha,\mu\>}
\end{equation}
whenever $\alpha\in\coroot\subset\lambdaz$. This binds $\omega$ uniquely
if $Z$ is cyclic, for if $\lambda\in\lambdaz$ is a generator of
$Z\cong\lambdaz/\coroot$, then
$\omega(\alpha+a\lambda,\beta+b\lambda)=
 (-1)^{\ell(\<\alpha,\beta\>+\<\lambda,a\beta+b\alpha\>)}$ for any
$\alpha,\beta\in\coroot$.
We therefore only need to investigate the case of
$G/Z=\SpinDD/\IZ_{2}\times\IZ_{2}$, $n\geq 2$. The main result
of this section  is that in this case $\ell$ is even and $\omega\equiv 1$.\\

We begin by computing a number of commutators related to
\eqref{eq:yet another}. Recall that the roots of $\SpinDD$
are the vectors $\pm\theta_{i}\pm\theta_{j}$, $1\leq i\neq j\leq 2n$,
where $\{\theta_{i}\}$ is the canonical basis of $\IR^{2n}$. The simple
roots are $\alpha_{i}=\theta_{i}-\theta_{i+1}$, $i=1\ldots 2n-1$ and
$\alpha_{2n}=\theta_{2n-1}+\theta_{2n}$ and the highest root is
$\theta=
 \theta_{1}+\theta_{2}=
 \alpha_{1}+2(\alpha_{2}+\cdots+\alpha_{2n-2})+\alpha_{2n-1}+\alpha_{2n}$
so that the minimal dominant coweights are $\cow{1}=\theta_{1}$,
$\cow{2n-1}=\half{1}(\theta_{1}+\cdots+\theta_{2n-1}-\theta_{2n})$
and $\cow{2n}=\half{1}(\theta_{1}+\cdots+\theta_{2n})$.
Fix, for any positive root $\alpha$, a basis 
$e_{\alpha},f_{\alpha},h_{\alpha}=\alpha^{\vee}$ of the corresponding
$\slc$--subalgebra of $\soDDc$ such that $e_{\alpha}^{*}=f_{\alpha}$
where $^{*}$ is the canonical anti--linear anti--involution acting
as -1 on $\soDD$. Then,

\begin{lemma}\label{le:weyl commutator}
Let $w_{j}$ be the Weyl group elements corresponding to the minimal
dominant coweights $\cow{j}$, $j\in\{1,2n-1,2n\}$ by proposition
\ref{centre weyl}. Then, the following elements may be taken
as representatives of $w_{j}$ in $\SpinDD$
\begin{align}
\wl{1}&=
\exp_{\SpinDD}\left(
\half{\pi}(e_{\theta_{1}+\theta_{2n}}-f_{\theta_{1}+\theta_{2n}})\right)
\exp_{\SpinDD}\left(
\half{\pi}(e_{\theta_{1}-\theta_{2n}}-f_{\theta_{1}-\theta_{2n}})\right)
\label{eq:w one}\\
\wl{2n-1}&=
\prod_{i=2}^{n}\exp_{\SpinDD}\left(
\half{\pi}(e_{\theta_{i}+\theta_{2n-i+1}}-f_{\theta_{i}+\theta_{2n-i+1}})\right)
\exp_{\SpinDD}\left(
\half{\pi}(e_{\theta_{1}-\theta_{2n}}-f_{\theta_{1}-\theta_{2n}})\right)\\
\wl{2n}&=
\prod_{i=1}^{n}\exp_{\SpinDD}\left(
\half{\pi}(e_{\theta_{i}+\theta_{2n-i+1}}-f_{\theta_{i}+\theta_{2n-i+1}})\right)
\label{eq:w twon}
\end{align}
Moreover, the group commutators
$[\wl{j},\wl{k}]=\comm{\wl{j}}{\wl{k}}$ are all equal to one.
\end{lemma}
\proof
It is readily seen that the $w_{j}$ may be expressed as
\begin{xalignat}{3}
w_{1}&=
\sigma_{\theta_{1}+\theta_{2n}}
\sigma_{\theta_{1}-\theta_{2n}}&
w_{2n-1}&=
\prod_{i=2}^{n}\sigma_{\theta_{i}+\theta_{2n-i+1}}
\sigma_{\theta_{1}-\theta_{2n}}&
w_{2n}&=
\prod_{i=1}^{n}\sigma_{\theta_{i}+\theta_{2n-i+1}}
\end{xalignat}
where $\sigma_{\alpha}$ is the orthogonal reflection corresponding
to the root $\alpha$ and acts on $\tc$ as
\begin{equation}\label{eq:reflection}
\sigma_{\alpha}(h)=h-\<h,\alpha\>\alpha^{\vee}
\end{equation}
Each $\sigma_{\alpha}$ may be lifted in $\SpinDD$ to
$\zigma{\alpha}=
\exp_{\SpinDD}\left(\half{\pi}(e_{\alpha}-f_{\alpha})\right)$.
Indeed, a power series expansion shows that $\Ad(\zigma{\alpha})$
leaves $\tc$ invariant and coincides with the right hand--side
of \eqref{eq:reflection}. Thus, \eqref{eq:w one}--\eqref{eq:w twon}
hold. The last claim is a consequence of the fact that the lifts of
$w_{j}$ and $w_{k}$ only involve roots $\alpha_{j_{p}}$, $\beta_{k_{q}}$
such that $\alpha_{j_{p}}\pm\beta_{k_{q}}$ is either zero or not a root.
Thus,
$[e_{\alpha_{j_{p}}}-f_{\alpha_{j_{p}}},
  e_{ \beta_{k_{q}}}-f_{ \beta_{k_{q}}}]=0$ and
$\zigma{\alpha_{j_{p}}}$ and $\zigma{\beta_{k_{q}}}$ commute
\halmos

\begin{lemma}\label{le:aut commutator}
Let $(\pi,\H)$ be an irreducible positive energy representation
of $L\SpinDD$ the isomorphism class of which is invariant under
$L_{\twotwo}\SpinDD$ and denote by the same symbol its unique
extension to the latter group. Then, if $\wl{j}$ are as in lemma
\ref{le:weyl commutator} and
$\zeta_{j}(\theta)=exp_{T}(-i\theta\cow{j})$ are the discontinuous
loops corresponding to the minimal dominant coweights $\cow{j}$
\begin{equation}
\pi(\zeta_{j}\wl{j})\pi(\zeta_{k}\wl{k})
\pi(\zeta_{j}\wl{j})^{*}\pi(\zeta_{k}\wl{k})^{*}=1
\end{equation}
\end{lemma}
\proof
Let $\cow{i}$ be the minimal dominant coweight with corresponding
Weyl group element $w_{i}$ such that $w_{j}w_{k}=w_{i}$. By
\eqref{compo}, $w_{j}\cow{k}=\cow{i}-\cow{j}$ so that
\begin{equation}\label{eq:tic}
 \wl{j}\zeta_{k}(\theta)\wl{j}^{-1}=
 \exp_{T}(-i\theta w_{j}\cow{k})=
 \zeta_{i}\zeta_{j}^{-1}(\theta)
\end{equation}
and similarly
\begin{equation}\label{eq:tac}
\wl{k}\zeta_{j}\wl{k}^{-1}=\zeta_{i}\zeta_{k}^{-1}
\end{equation}
Thus,
\begin{equation}
[\zeta_{j}\wl{j},\zeta_{k}\wl{k}]=
 \zeta_{j}\wl{j}\zeta_{k}\wl{j}^{-1}[\wl{j},\wl{k}]
 \wl{k}\zeta_{j}^{-1}\wl{k}^{-1}\zeta_{k}^{-1}=
 \zeta_{i}[\wl{j},\wl{k}]\zeta_{i}^{-1}
\end{equation}
which, by lemma \ref{le:weyl commutator}, equals 1. It follows that
the group commutator $[\pi(\zeta_{j}\wl{j}),\pi(\zeta_{k}\wl{k})]$ acts
as a scalar $\chi$ on $\H$. To evaluate $\chi$, recall that by corollary
\ref{automorphism}, conjugation by $\zeta_{j}\wl{j}$ and $\zeta_{k}\wl{k}$
induces an automorphism of $\lpol\soDD$ preserving its triangular
decomposition so that the unitaries $\pi(\zeta_{j}\wl{j}),
\pi(\zeta_{k}\wl{k})$ leave the highest weight vector in $\H$ invariant
and $\chi=1$ \halmos

\begin{theorem}\label{th:mackey obs}
Let $\pi$ be an irreducible positive energy representation of $L\SpinDD$
the isomorphism class of which is invariant under $L_{\twotwo}\SpinDD$
and extend it uniquely to a projective representation of the latter group.
Then, the corresponding central extension of $L_{\twotwo}\SpinDD$ has even
level and trivial commutator map.
\end{theorem}
\proof
As readily checked from equation \eqref{eq:n}, the geometric action of
$Z(\SpinDD)$ on the level $\ell$ alcove has fixed points only if $\ell$
is even. Thus from \eqref{triviality}, we get $\omega(\alpha,\mu)=1$
whenever $\alpha\in\coroot$ and $\omega$ is the pull--back of one of
the two skew--symmetric forms on
$\coweight/\coroot\cong\IZ_{2}\times\IZ_{2}$ so that $\omega\equiv 1$
iff $\omega(\cow{j},\cow{k})=1$ for a pair of distinct minimal dominant
coweights $\cow{j},\cow{k}$.\\

By lemma \ref{le:weyl commutator} and the fact that $\pi$ lifts to a
unitary representation of $\SpinDD$
\begin{equation}\label{eq:rewrite}
\begin{split}
[\pi(\zeta_{j}\wl{j}),\pi(\zeta_{k}\wl{k})]
&=
 \pi(\zeta_{j})\pi(\zeta_{k})[\pi(\zeta_{k})^{*},\pi(\wl{j})]
[\pi(\wl{j}),\pi(\wl{k})]
[\pi(\wl{k}),\pi(\zeta_{j})^{*}]\pi(\zeta_{j})^{*}\pi(\zeta_{k})^{*}\\
&=
 \pi(\zeta_{j})\pi(\zeta_{k})[\pi(\zeta_{k})^{*},\pi(\wl{j})]
[\pi(\wl{k}),\pi(\zeta_{j})^{*}]\pi(\zeta_{j})^{*}\pi(\zeta_{k})^{*}
\end{split}
\end{equation}
By \eqref{eq:tic},
$[\zeta_{k}^{-1},\wl{j}]=\zeta_{k}^{-1}\zeta_{i}\zeta_{j}^{-1}$ so that
$[\pi(\zeta_{k})^{*},\pi(\wl{j})]$ is proportional to
$\pi(\zeta_{k})^{*}\pi(\zeta_{i})\pi(\zeta_{j})^{*}$. Similarly, by
\eqref{eq:tac}, $[\pi(\wl{k}),\pi(\zeta_{j})^{*}]$ is proportional
to $\pi(\zeta_{k})\pi(\zeta_{i})^{*}\pi(\zeta_{j})$ so that, by
\eqref{eq:yet another}, \eqref{eq:rewrite} is equal to
\begin{equation}
\omega(\cow{j},\cow{k})
[\pi(\zeta_{k})^{*},\pi(\wl{j})]
[\pi(\wl{k}),\pi(\zeta_{j})^{*}]
\end{equation}
and, by lemma \ref{le:aut commutator},
\begin{equation}\label{eq:obstruction}
\omega(\cow{j},\cow{k})=
[\pi(\zeta_{j})^{*},\pi(\wl{k})]
[\pi(\wl{j}),\pi(\zeta_{k})^{*}]
\end{equation}
We shall now show that the right hand--side of \eqref{eq:obstruction}
is equal to one.\\

Set $j=1$ and $k=2n$ so that the corresponding coweights
are $\cow{1}=\theta_{1}$ and
$\cow{2n}=\half{1}(\theta_{1}+\cdots+\theta_{n})$. As previously
noted, if $x_{\beta}\in\soDDc$ is a root vector corresponding to
$\beta$ and $\zeta_{\lambda}(\theta)=\exp_{T}(-i\theta\lambda)$,
then
\begin{equation}
\zeta_{\lambda}^{-1}x_{\beta}(n)\zeta_{\lambda}=
 x_{\beta}(n+\<\lambda,\beta\>)
\end{equation}
and therefore
\begin{equation}
\begin{split}
\zeta_{1}^{-1}\wl{2n}\zeta_{1}
&=
\exp_{L\SpinDD}\left( \half{\pi}
(e_{\theta_{1}+\theta_{2n}}(1)-f_{\theta_{1}+\theta_{2n}}(-1))\right)\\
&\prod_{i=2}^{n}
\exp_{ \SpinDD}\left( \half{\pi}
(e_{\theta_{i}+\theta_{2n-i+1}}-f_{\theta_{i}+\theta_{2n-i+1}})
\right)
\end{split}
\end{equation}
\begin{equation}
\begin{split}
\zeta_{2n}^{-1}\wl{1}^{-1}\zeta_{2n}
=&
\exp_{ \SpinDD}\left(-\half{\pi}
(e_{\theta_{1}-\theta_{2n}}   -f_{\theta_{1}-\theta_{2n}}    )\right)\\
 &
\exp_{L\SpinDD}\left(-\half{\pi}
(e_{\theta_{1}+\theta_{2n}}(1)-f_{\theta_{1}+\theta_{2n}}(-1))\right)
\end{split}
\end{equation}

Denoting by $\rho$ the infinitesimal action of $\lpol\soDDc$ on
$\hfin$ corresponding to $\pi$ via theorem \ref{th:core}, we
therefore get by proposition \ref{conjugated action}

\begin{align}
[\pi(\zeta_{1})^{*},\pi(\wl{2n})]&=
\exp\left( \half{\pi}\rho
(e_{\theta_{1}+\theta_{2n}}(1)-f_{\theta_{1}+\theta_{2n}}(-1))\right)
\exp\left(-\half{\pi}\rho
(e_{\theta_{1}+\theta_{2n}}   -f_{\theta_{1}+\theta_{2n}}    )\right)\\
[\pi(\wl{1}),\pi(\zeta_{2n})]&=
\exp\left( \half{\pi}\rho
(e_{\theta_{1}+\theta_{2n}}   -f_{\theta_{1}+\theta_{2n}}    )\right)
\exp\left(-\half{\pi}\rho
(e_{\theta_{1}+\theta_{2n}}(1)-f_{\theta_{1}+\theta_{2n}}(-1))\right)
\end{align}
where the exponentials are given by the spectral theorem. Thus,
\begin{equation}
[\pi(\zeta_{1})^{*},\pi(\wl{2n})][\pi(\wl{1}),\pi(\zeta_{2n})]=1
\end{equation}
as claimed \halmos

\ssection{Positive energy representations of $\lzg$}
\label{se:per of lzg}
%===================================================

Let $k$ be the order of the largest cyclic subgroup of $Z$ and
consider the action of the universal $k$--coverings of $\diff$
and $\rot$ on $\lzg$ by reparametrisation, as in \S \ref{ss:action of diff}.
We denote these coverings by $\diffk$ and $\rotk$ respectively.
Define a positive energy representation of $\lzg$ to be a strongly
continuous homomorphism
\begin{equation}
\pi:\lzg\longrightarrow PU(\H)
\end{equation}
extending to $\lzg\rtimes\rotk$ in such a way that $\rotk$ acts
by non--negative characters only and with finite--dimensional
eigenspaces. 

\begin{theorem}\label{th:classification of per}
An irreducible positive energy representation $(\pi,\H)$
of $\lzg$ yields by restriction
\begin{enumerate}
\item A positive energy representation of $LG$ of level
\begin{equation}\label{eq:constr one}
\ell\in\ell_{f}\cdot\IN
\end{equation}
where $\ell_{f}\in\{1,2\}$ is the fundamental level of $G/Z$.
\item A projective representation of $\lambdaz\cong\Hom(\T,T/Z)$,
the commutator map of which, defined by
$\omega(\lambda,\mu)=
 \pi(\zeta_{\lambda})\pi(\zeta_{\mu})
 \pi(\zeta_{\lambda})^{*}\pi(\zeta_{\mu})^{*}$,
satisfies
\begin{equation}\label{eq:constr two}
\omega(\alpha,\mu)=(-1)^{\ell\<\alpha,\mu\>}
\end{equation}
whenever $\alpha$ lies in the coroot lattice $\coroot\cong\Hom(\T,T)$.
\end{enumerate}
As an $LG$--module,
\begin{equation}\label{eq:orbit}
\H=\bigoplus_{\mu\in Z\lambda}\H_{\mu}\otimes\IC^{m_{\lambda}}
\end{equation}
where $\lambda$ lies in the level $\ell$ alcove of $G$, $Z\lambda$
is its orbit under the action of $Z\subseteq Z(G)$ defined by
proposition \ref{geometric action} and $\H_{\mu}$ is the irreducible
level $\ell$ positive energy representation of $LG$ with highest weight
$\mu$. Moreover, $m_{\lambda}=1$ unless $G/Z=\PSO_{4n}$,
$Z\lambda=\{\lambda\}$, $\ell$ is even and $\omega$ is the pull--back
of the non--trivial, skew--symmetric form on $Z\cong\IZ_{2}\times\IZ_{2}$,
in which case $m_{\lambda}=2$. The triple $(\ell,\omega,Z\lambda)$
classifies $\H$ uniquely and for any $(\ell,\omega,Z\lambda)$ satisfying
\eqref{eq:constr one} and \eqref{eq:constr two}, there exists
an irreducible positive energy representation of $\lzg$ realising
it. Lastly, the action of $\rotk$ on $\H$ extends uniquely to
a projective unitary representation $\rho$ of $\diffk$ satisfying
\begin{equation}\label{eq:intertwine}
\rho(\phi)\pi(\zeta)\rho(\phi)^{*}=\pi(\zeta\circ\phi^{-1})
\end{equation}
which coincides with the Segal--Sugawara representation obtained
by regarding $\H$ as a positive energy representation of $LG$.
\end{theorem}
\proof
We shall repeatedly, and without further mention, use the following
fact. Let $Y\subseteq Z$ be a subgroup and $(\rho,\K)$ a positive
energy representation of $\lyg$. $\rho$ lifts to a unitary representation
of the continuous central extension $\pbackdue$ of $\lyg$ obtained
by pulling back
\begin{equation}
1\rightarrow\T\rightarrow U(\K)\xrightarrow{p} PU(\K)\rightarrow 1
\end{equation}
to $\lzg$. Explicitly,
\begin{equation}
\pbackdue=\{(\zeta,V)\in\lyg\times U(\K)|\rho(\zeta)=p(V)\}
\end{equation}
acts on $\K$ by $(\zeta,V)\xi=V\xi$. Let $h$ be the level of $\K$
as a positive energy representation of $LG$. Then, the restriction
of $\pbackdue$ to $LG$ is smooth and of level $h$ \cite[\S II.2.4]{TL}
so that $\pbackdue$ is a smooth central extension of $\lyg$ and
therefore falls within the classification of section
\ref{se:extensions of lzg}. Set now $Y=Z$, $\rho=\pi$ and $\K=\H$.
Then, the level $\ell$ and commutator map $\omega$ of $\pback$,
which is readily seen to be
\begin{equation}
 \omega(\lambda,\mu)=
 \pi(\zeta_{\lambda})\pi(\zeta_{\mu})
 \pi(\zeta_{\lambda})^{*}\pi(\zeta_{\mu})^{*}
\end{equation}
are bound by theorem \ref{th:constraint} and therefore satisfy
\eqref{eq:constr one}--\eqref{eq:constr two}.\\

As an $LG\rtimes\rotk$--module, $\H$ decomposes as
\begin{equation}\label{eq:decomp one}
\H=\bigoplus_{\mu}\H(\mu)
\end{equation}
where $\mu$ spans the level $\ell$ alcove $\al$ of $G$ and $\H(\mu)$
is the isotypical summand of $\H$ corresponding to the irreducible,
level $\ell$ positive energy representation $\H_{\mu}$ of $LG$ with
highest weight $\mu$. Evidently, for any $\zeta\in\lzg$,
$\pi(\zeta)\H(\mu)=\H(\zeta\mu)$ where the notation refers to
the abstract action of $Z\cong\lzg/LG$ on $\al$ given by propositions
\ref{pr:stability} and \ref{conjugated action} which, by theorem
\ref{th:Z on per}, coincides with the geometric one given by proposition
\ref{geometric action}. Since $\H$ is irreducible, \eqref{eq:decomp one}
reduces to
\begin{equation}
\H=\bigoplus_{\mu\in Z\lambda}\H(\mu)
\end{equation}
for some $\lambda\in\al$ and the triple $(\ell,\omega,Z\lambda)$
is an invariant of $\H$.\\

To proceed, it will be more convenient to consider unitary
representations rather than projective ones. For any subgroup
$Y\subseteq Z$ and central extension $\wt\lyg$ of $\lyg$ with
level $h$ and commutator map $\kappa$, there is a bijective
correspondence between positive energy representations of $\lyg$
corresponding to the pair $(h,\kappa)$ and unitary representations
of $\wt\lyg\rtimes\rotk$ such that $\rotk$ acts by non--negative
characters only and with finite--dimensional eigenspaces and
the central subgroup $\T\subset\wt\lyg$ acts as multiplication
by the character $z\rightarrow z$, provided that representations
differing by a character of $\wt\lyg$ are identified \footnote
{Since $LG$ is perfect \cite[prop. 3.4.1]{PS}, such characters
factor through the group of components $Y$ of $\lyg$.}. We call
these positive energy representations of $\wt\lyg$ and will work
with them from now on.\\

Fix now $\ell,\omega$ satisfying
\eqref{eq:constr one}--\eqref{eq:constr two} and denote by
$\wt\lzg$ the central extension of $\lzg$ with level $\ell$
and commutator map $\omega$, the existence of which is
guaranteed by proposition \ref{pr:ext of lzg}. Recall moreover
that by proposition \ref{pr:action of diff}, the action of
$\diffk$ on $\lzg$ lifts to $\wt\lzg$.
For a given orbit $Z\lambda\subseteq\al$ with isotropy subgroup
$Y\subseteq Z$, denote by $\wt\lyg$ and $\wt{LG}$ the restrictions
of $\wt\lzg$ to $\lyg$ and $LG$ respectively. Then, by Mackey's
theory \cite[thm 3.11]{Ma2}, and in view of the fact that positive
energy representations of $LG$ are invariant under conjugation by
$\lzg$, the map
\begin{equation}
i:\K\longrightarrow
  \ind{\wt\lyg\rtimes\rotk}{\wt\lzg\rtimes\rotk}\K
\end{equation}
gives a bijection between the irreducible positive energy
representations $(\rho,\K)$ of $\wt\lyg$ the restriction to
$\wt{LG}$ of which is isotypical of type $\H_{\lambda}$, and
the irreducible positive energy representations of $\wt\lzg$
with highest weight orbit $Z\lambda$ \footnote
{The induction functor is well--defined in the present
context since $\wt\lyg\rtimes\rotk\subseteq\wt\lzg\rtimes\rotk$
is of finite index, and satisfies the usual properties
of its finite--dimensional counterpart which are necessary
to prove Mackey's theorem.
Moreover, an elementary application of the induction--restriction
theorem shows that $\T\subset\wt\lzg$ acts by the required
character on $i(\K)$.}.
Moreover, since any character $\chi$ of $Y$ extends to $Z$ and
$\ind{}{}(\K\otimes\chi)=\ind{}{}(\K)\otimes\chi$, $\K$ and $\K'$
differ by a character iff $i(\K)$ and $i(\K')$ do.\\

Notice also that $i(\K)$ admits an intertwining action of $\diffk$ if
$\K$ does. Indeed, let $R$ be a projective representation of $\diffk$
on $\K$ satisfying
\begin{equation}
R(\phi)\rho(\wt\zeta)R(\phi)^{*}=\rho(\wt\zeta\circ\phi^{-1})
\end{equation}
for any $\wt\zeta\in\wt\lyg$, and lift it to a unitary
representation of the corresponding central extension
$\wt\diffk=R^{*}U(\K)$ of $\diffk$. Then, by induction--restriction,
\begin{equation}
i(\K)=
\ind{\wt\lyg\rtimes\rotk}{\wt\lzg\rtimes\rotk}\K\cong
\ind{\wt\lyg\rtimes\wt\diffk}{\wt\lzg\rtimes\wt\diffk}\K
\end{equation}
as $\wt\lzg\rtimes\rotk$--modules, and $i(\K)$ admits an intertwining
action of $\diffk$.\\

The relevant representations of $\wt\lyg$ are obtained in the
following way \cite[\S 3.10]{Ma2}. Extend the unitary action
$\pil$ of $\wt{LG}\rtimes\rotk$ on $\H_{\lambda}$ to a projective
one of $\wt\lyg\rtimes\rotk$ satisfying, and uniquely determined
by
\begin{equation}\label{eq:meta}
\pil(\wt\zeta)\pil(\wt\gamma)\pil(\wt\zeta)^{*}=
\pil(\wt\zeta\wt\gamma{\wt\zeta}^{-1})
\end{equation}
for any $\wt\zeta\in\wt\lyg$ and $\wt\gamma\in\wt{LG}$.
The corresponding central extension $\pil^{*}U(\H_{\lambda})$
of $\wt\lyg$ determines one of $Y$ by
\begin{equation}
\wt Y=
\pil^{*}U(\H_{\lambda})/
(\wt\gamma,\pil(\wt\gamma))_{\wt\gamma\in\wt{LG}}
\end{equation}
Any irreducible unitary representation $\rho$ of $\wt Y$ such
that its central subgroup $\T$ acts by $z\rightarrow z^{-1}$
yields one of $\pil^{*}U(\H_{\lambda})$, namely $\pil\otimes\rho$,
where $\T$ acts trivially, and therefore one of 
$\wt\lzg\cong\pil^{*}(U(\H_{\lambda}))/\T$ and the representations
of $\wt\lyg$ in question are exactly of this form.\\

Moreover, they admit an intertwining action of $\diffk$. Indeed,
if $R$ is the Segal--Sugawara representation of $\diff$ on
$\H_{\lambda}$ intertwining $\wt{LG}$ \cite[prop. 13.4.2]{PS},
then, for any $\wt\zeta\in\wt\lyg$ and $\phi\in\diffk$
\begin{equation}
R(\phi)\pil(\wt\zeta)R(\phi)^{*}=\pil(\wt\zeta\circ\phi^{-1})
\end{equation}
projectively, since both sides have the same commutation relations
with $\wt{LG}$. Thus,
\begin{equation}
(R(\phi)\otimes 1)
\pil\otimes\rho(\wt\zeta)
(R(\phi)\otimes 1)^{*}=
\kappa(\phi,\wt\zeta)
\pil\otimes\rho(\wt\zeta\circ\phi^{-1})
\end{equation}
for some $\kappa(\phi,\wt\zeta)\in\T$ which is multiplicative
in each variable. Since $\diffk$ is perfect \cite{He,Ep},
$\kappa\equiv 1$ and $\pil\otimes\rho$ admits an intertwining
action of $\diffk$.\\

To determine $\wt Y$, notice that the projective unitaries
$\pil(\wt\zeta)$, $\wt\zeta\in\wt\lyg$ defined by \eqref{eq:meta}
only depend on the image of $\wt\zeta$ in $\lyg$ and therefore
give rise to a central extension of $\lyg$. Since this extension
has level $\ell$, it differs from $\wt\lyg$ by the pull--back
of an extension of $Y$ which is readily seen to be $\wt Y$.
We must now distinguish two cases. If $\wt Y$ splits, which
is so if $Y$ is cyclic or, by theorem \ref{th:mackey obs}, if
$Y=Z=Z(\Spin_{4n})$ and $\omega\equiv 1$, the relevant
representations of $\wt Y$ correspond to the characters $\chi$
of $Y$ and the irreducible, positive energy representations of
$\wt\lzg$ with highest weight orbit $Z\lambda$ are of the form
\begin{equation}
\ind{\wt\lyg\rtimes\rotk}{\wt\lzg\rtimes\rotk}(\pil\otimes\chi)
\end{equation}
so that their restriction to $LG$ only involves isotypical
summands of multiplicity one.
If, on the other hand, $\wt Y$ doesn't split, then
$Y=Z=Z(\Spin_{4n})\cong\IZ_{2}\times\IZ_{2}$, $\ell$ is even
and $\omega$ is the pull--back of the non--trivial,
skew--symmetric form on $Y$. In this case, the relevant
representation is the Heisenberg representation of $\wt Y$ on
$\K\cong\IC^{2}$ and there exists a unique irreducible, positive
energy representation of $\lzg$ with highest weight orbit
$Z\lambda=\{\lambda\}$ the restriction to $LG$ of which
is isomorphic to $\H_{\lambda}\otimes\IC^{2}$.\\

This completes the classification of irreducible positive energy
representations of $\lzg$ and shows that any such $(\pi,\H)$
admits an intertwining action $\rho$ of $\diffk$. Let now $\rho_{0}$
be the Segal--Sugawara representation of $\diffm$ on $\H$, where
$m$ is some positive integer which we may take to be a multiple
of $k$
\footnote{The Segal--Sugara representation factors through
$\diff$ only on irreducible positive energy representations
of $LG$ and through a finite order covering of $\diff$ in
general.}. To see that $\rho$ coincides with $\rho_{0}$, notice that
for any $\phi\in\diffm$, the operator $\rho(\phi)\rho_{0}(\phi)^{*}$
commutes with $LG$. Thus, if
\begin{equation}
\H=\bigoplus_{\mu\in Z\lambda}\H_{\mu}\otimes\IC^{m_{\lambda}}
\end{equation}
is the decomposition of $\H$ as an $LG$--module,
then $\rho_{0}(\phi)=\bigoplus_{\mu}\rho_{0}^{\mu}(\phi)\otimes 1$
where $\rho_{0}^{\mu}(\phi)$ gives the Segal--Sugawara representation
of $\diff$ on $\H_{\mu}$ and therefore
\begin{equation}
\rho(\phi)=\bigoplus_{\mu}\rho_{0}^{\mu}(\phi)\otimes T_{\mu}(\phi)
\end{equation}
for some unitary operators $\T_{\mu}(\phi)$. Since both $\rho$ and
the $\rho_{0}^{\mu}$ are projective representations of $\diffm$, the
same is true of the $T_{\mu}$ which are therefore trivial since
$\diffm$ admits no finite--dimensional projective representations
\cite[prop. 3.3.2]{PS} \halmos

\ssection{Positive energy representations of $L(G/Z)$}
\label{se:factoring}
%====================================================

We now determine those positive energy representations of $\lzg$ which
factor through $L(G/Z)\cong\lzg/Z$.

\begin{lemma}\label{le:char id}
Let $(\pi,\H)$ be an irreducible positive energy representation of $LG$
with highest weight $\lambda$ and consider its unique lift to a unitary
representation of $G$. Then, $Z(G)$ acts on $\H$ as multiplication by
the character $\chi(\exp_{T}(h))=e^{\<\lambda,h\>}$.
\end{lemma}
\proof
For any $z\in Z(G)$ and $\gamma\in LG$ we have
$\pi(\gamma)\pi(z)\pi(\gamma)^{*}\pi(z)^{*}=\kappa(\gamma,z)$ where
$\kappa(\cdot,\cdot)\in T$ is independent of the particular choice
of lifts. $\kappa$ is continuous and multiplicative in both variables
and therefore defines a continuous map $LG\rightarrow\widehat{Z(G)}$
which, by connectedness of $LG$, is trivial. Thus, by Shur's lemma,
$Z(G)$ as multiplication by a character which is easily computed since
if $\Omega\in\H$ is the highest weight vector in $\H$, then
$\pi(\exp_{T}(h))\Omega=e^{\<\lambda,h\>}\Omega$ \halmos\\

Let now $Z\cong\lambdaz/\coroot$ be a subgroup of $Z(G)$. By using the
basic inner product $\<\cdot,\cdot\>$ on $i\t$, the dual group $\wh{Z}$
may be identified with $\weight/(\lambdaz)^{*}$ where
$\root\subseteq(\lambdaz)^{*}\subseteq\weight$ is the dual lattice of
$\lambdaz$.

\begin{lemma}
Let $\H$ be an irreducible positive energy representation of $\lzg$ of
level $\ell$ and highest weight orbit $Z\lambda$. Then, the characters
of $Z\subset\lzg$ corresponding to $\H$ are the classes of
$\lambda+\ell\cow{i}$ mod $(\lambdaz)^{*}$ where $\cow{i}$
are the minimal dominant coweights corresponding to $Z$.
\end{lemma}
\proof
When restricted to $LG$, $\H$ decomposes as a direct sum of positive
energy representations of $LG$ the highest weights of which lie on
the orbit $Z\lambda$. By lemma \ref{le:char id} and proposition
\ref{geometric action}, these give rise to the characters
$\ell\cow{i}+w_{i}\lambda$ mod $(\lambdaz)^{*}$ of $Z$ where $\cow{i}$
are the minimal dominant coweights corresponding to $Z$. Since $W$
preserves $\root$, and {\it a fortiori} $(\lambdaz)^{*}$--cosets in
$\weight$ however, we get
$\ell\cow{i}+w_{i}\lambda=\ell\cow{i}+\lambda$ mod $(\lambdaz)^{*}$
\halmos

\begin{corollary}
An irreducible positive energy representation $\pi$ of $\lzg$ factors
through $L(G/Z)$ if, and only if its level is a a multiple of the basic
level $\ellb$ of $G/Z$.
\end{corollary}
\proof
$\pi$ factors through $L(G/Z)$ iff $Z$ acts by the same character on
each of its irreducible $LG$--submodules. By definition, $\ellb$ is
the smallest integer $\ell$ such that $\ell\<\cdot,\cdot\>$ is integral
on $\lambdaz$ and therefore such that $\ell\cow{i}\in(\lambdaz)^{*}$
for any fundamental coweight $\cow{i}$ corresponding to $Z$ \halmos


\begin{thebibliography}{ZZZ}
%======================

% \bibitem[Ad]{Ad} J. F. Adams, {\it Lectures on Lie Groups.}
% Univ. of Chicago Press, 1969.

\bibitem[Bou]{Bou} N. Bourbaki, {\it Groupes et Alg\`ebres de Lie,
Chapitres 4,5 et 6}. Masson, Paris 1981.

\bibitem[Ep]{Ep} D. B. A. Epstein, {\it Commutators of
$C^{\infty}$--Diffeomorphisms. Appendix to ``A Curious Remark
Concerning the Geometric Transfer Map'' by John Mather}, Comment.
Math. Helv. 59 (1984), 111-122.

\bibitem[FGK1]{FGK1} G. Felder, K. Gaw\c{e}dzki, A. Kupiainen, {\it The
Spectrum of Wess--Zumino--Witten Models}, Nucl. Phys. B 299 (1988),
355-66.

\bibitem[FGK2]{FGK2} G. Felder, K. Gaw\c{e}dzki, A. Kupiainen, {\it
Spectra of Wess--Zumino--Witten Models with Arbitrary Simple Groups},
Comm. Math. Phys. 117 (1988), 127--58.

\bibitem[Ga]{Ga} M. R. Gaberdiel, {\it WZW Models of General Simple Groups},
Nucl. Phys. B 460 (1996), 181--202.

\bibitem[GW]{GW} D. Gepner, E. Witten, {\it String Theory on Group
Manifolds}, Nucl. Phys. B 278 (1986), 493--549.

\bibitem[GO]{GO} P. Goddard, D. Olive, {\it The Magnetic Charges of Stable
Self--Dual Monopoles}, Nucl. Phys. B191 (1981), 528--548.
303-414.

%\bibitem[GO]{GO} P. Goddard, D. Olive, {\it Kac-Moody and Virasoro Algebras in
%Relation to Quantum Physics}, International Journal of Modern Physics A 1(1986)
%303-414.

\bibitem[GoWa]{GoWa} R. Goodman, N. R. Wallach, {\it Structure and Unitary
Cocycle Representations of Loop Groups and the Group of Diffeomorphisms
of the Circle}, J. Reine Angew. Math. 347 (1984), 69--133.

\bibitem[He]{He} M.-R. Herman, {\it Simplicit\'e du Groupe des
Diff\'eomorphismes de Classe $C^{\infty}$, Isotopes \`a l'Identit\'e,
du Tore de Dimension $n$}, C. R. Acad. Sci. Paris, Ser. A, 273 (1971),
232--4.

\bibitem[Hu]{Hu} J. E. Humphreys, {\it Introduction to Lie Algebras and 
Representation Theory.} Graduate Texts in Mathematics, 
Springer-Verlag 1972.

\bibitem[Ka]{Ka} V. Kac, {\it Infinite Dimensional Lie Algebras}, 3rd. edition,
Cambridge University Press, 1994.

\bibitem[Ma1]{Ma1} G. W. Mackey, {\it Unitary Representations of Group
Extensions I}, Acta. Math. 99 (1958) 265--311.

\bibitem[Ma2]{Ma2} G. W. Mackey, {\it The Theory of Unitary Group
Representations}. University of Chicago Press, 1976.

\bibitem[OT]{OT} D. Olive, N. Turok, {\it The Symmetries of Dynkin Diagrams
and the Reduction of Toda Field Equations}, Nucl. Phys. B{\bf 215} (1983),
470-94.

\bibitem[PS]{PS} A. Pressley, G. B. Segal, {\it Loop Groups}. Oxford 
University Press 1986.

\bibitem[Wa]{Wa} A. J. Wassermann, {\it Conformal Field Theory and Operator 
Algebras II : Fusion for Von Neumann Algebras and Loop groups}, preprint.

\bibitem[TL]{TL} V. Toledano Laredo, {\it Fusion of Positive Energy
Representations of LSpin$_{2n}$}, Ph.D. dissertation, University of
Cambridge, 1997.

\bibitem[TL2]{TL2} V. Toledano Laredo, {\it Integrating Unitary
Representations of Infinite-Dimensional Lie Groups}. J. Funct.
Anal. {\bf 161} (1999), 478--508.

\end{thebibliography}
\end{document}